\documentclass[sts]{imsart}
\RequirePackage{amsthm,amsmath,amsfonts,amssymb}
\RequirePackage[numbers]{natbib}
\RequirePackage[colorlinks,citecolor=blue,urlcolor=blue]{hyperref}
\RequirePackage{graphicx}

\usepackage[utf8]{inputenc}
\usepackage{adjustbox}

\usepackage{latexsym}
\usepackage{enumerate}
\usepackage{wasysym}
\usepackage{lipsum}
\usepackage{tcolorbox}

\usepackage[graphicx]{realboxes}

\usepackage{float}
\usepackage{lscape}

\usepackage{booktabs,color,epsfig,graphicx}
\usepackage{subfigure}
\usepackage{soul}

\startlocaldefs 

\newtheorem{guideline}{Guideline}

\hypersetup{urlcolor=blue, citecolor=blue, colorlinks=true, linkcolor=blue} 

\DeclareMathOperator*{\argmin}{arg\,min}

\def\R{\mathbb{R}}
\def\Z{\mathbb{Z}}  

\def\b_eta{\boldsymbol{\eta}}

\def\calH{\cal H}  
\def\bZ{\boldsymbol{Z}} 
\def\bR{\boldsymbol{R}} 
\def\bx{\boldsymbol{x}} 
\def\by{\boldsymbol{y}} 

\def\bm{\boldsymbol}
\def\M{Mat\'ern {}} 
\endlocaldefs

\begin{document}

\begin{frontmatter}
\title{The Mat{\'e}rn Model: 
\\ A Journey through Statistics, Numerical Analysis and Machine Learning}
\runtitle{Mat{\'e}rn: a Journey}
\begin{aug}
\author[A]{\fnms{Emilio}~\snm{Porcu}\thanksref{t1}\ead[label=e1]{emilio.porcu@bku.ac.ae}},
\author[B]{{Moreno} {Bevilacqua}\ead[label=e2]{moreno.bevilacqua@uai.cl}},
\author[C]{{Robert} {Schaback}\ead[label=e3]%
{schaback@math.uni-goettingen.de}}
\and
\author[D]{{Chris J.} {Oates}\ead[label=e4]{chris.oates@ncl.ac.uk}} 
\thankstext{t1}{{Corresponding Author.}}

\address[A]{Emilio Porcu is Professor, Department of Mathematics,
Khalifa University, Abu Dhabi, UAE \printead{e1}.}

\address[B]{Moreno Bevilacqua is Professor, Department of Statistics, Universidad Adolfo Ibanez \printead{e2}.}

\address[C]{Robert Schaback is Professor, Department of Mathematics, University of G\"ottingen, Germany \printead[presep={\ }]{e3}.} 
\address[D]{Chris J. Oates is Professor, School of Mathematics, Statistics \& Physics, Newcastle University, UK \printead{e4}.}
\end{aug}

\begin{abstract}
The Mat{\'e}rn model has been a cornerstone of spatial statistics for more than half a century. More recently, the Mat{\'e}rn model has been central to disciplines as diverse as 
numerical analysis, approximation theory, computational statistics, machine learning, and probability theory. In this article we take a Mat{\'e}rn-based journey across these disciplines. 
First, we reflect on the importance of the Mat{\'e}rn model for estimation and prediction in spatial statistics, establishing also connections to other disciplines in which the Mat{\'e}rn model has been influential. 
Then, we position the Mat{\'e}rn model within the literature on big data and scalable computation: the SPDE approach, the Vecchia likelihood approximation, and recent applications in Bayesian computation are all discussed. 
Finally, we review recent devlopments, including flexible alternatives to the Mat{\'e}rn model, whose performance we compare in terms of estimation, prediction, screening effect, computation, and Sobolev regularity properties. 


\vspace{1cm}

\end{abstract}
\end{frontmatter}
\newpage


\section{Introduction} \label{sec1}


This paper serves two purposes: 
On the one hand, we provide a panoramic view, across several disciplines, of the Mat{\'e}rn model.
On the 
other 
hand, the paper provides constructive criticisms about the role of the Mat{\'e}rn model in several disciplines, while discussing alternative or more general models and their relevance to many aspects of statistical modeling, estimation, prediction, computational statistics, numerical analysis, and machine learning.  

A historical account of the Mat{\'e}rn model is provided by \citet{Guttorp}. 
The Mat{\'e}rn model -- also called the Mat{\'e}rn \emph{covariance function}, or the Mat{\'e}rn \emph{kernel}, depending on context -- is commonly attributed to \citet{Matern:1960}, but can be found under alternative names in different branches of the scientific literature. 
The use of the Mat{\'e}rn model 
is widespread, and it is impossible to cover all its diverse applications here; our review focuses on a selection of applications that are of especial interest and significance. 
Specifically, we aim to cover
\begin{enumerate}
\item 
estimation and prediction using the Mat{\'e}rn model in statistics, with emphasis on 
maximum likelihood estimation, 
Kriging 
prediction, 
and the associated screening effect;

\item applications of the Mat{\'e}rn model in
\begin{enumerate}
\item computational statistics, 
including the stochastic differential equation (SDE) and stochastic partial differential equation (SPDE) approaches, likelihood approximation, inference of partial differential equations (PDEs) and Stein's method; 
\item statistical modeling, 
including non-standard scenarios, for instance when isotropy and stationarity cannot be assumed,
or to model directions and curves;
\item approximation theory and numerical analysis, where the Mat{\'e}rn model is used to construct 
kernel-based interpolants; 
\item machine learning, 
where the Mat{\'e}rn model 
is central to the literature on Gaussian processes modelling; and
\item probability theory, where the Mat{\'e}rn model has inspired several contributions based on properties of the sample paths of associated stochastic processes, in concert with the solution of certain classes of stochastic differential equations;
\end{enumerate}
\item 
comparison with recent flexible alternatives  
to the Mat{\'e}rn model, with a focus on
\begin{enumerate}
\item enhanced models with interesting features, such as compact support or polynomial decay;
\item asymptotic estimation accuracy, misspecified prediction, and screening effects;
\item the implications of using certain classes of compactly supported kernels within approximation theory, computational statistics, and machine learning. 
\end{enumerate}
\end{enumerate}
This article is novel, in being the first to take a broad view of the scientific literature through the lens of the Mat{\'e}rn model.
In particular, we do \textit{not} attempt a review of covariance functions in general. 
Recent reviews provide a quite exhaustive panorama of covariance models, from space to space-time \citep{porcu202130}, to multivariate covariance functions \citep{Genton:Kleiber:2014}, and covariance-based modeling on spheres and manifolds \citep{PAF2016}.
In addition, while there are many fascinating applications of the Mat{\'e}rn model across the scientific landscape, we cannot hope to do justice to them all.
Our emphasis is therefore limited to methodological and theoretical issues which we hope are of relevance across a wide range of disciplines in which the Mat{\'e}rn model is used.



\subsection{Setting and Notation} \label{background}


Throughout, bold letters refer to vectors and matrices, and the transpose operator is denoted $\top$. 
Let $d \in \mathbb{N}$ and let $Z=\{{Z}(\bm{x}), \; \bm{x} \in \R^d \}$ be a real-valued Gaussian random field, having zero mean and and \emph{covariance function} ${K}: \R^d \times \R^d \to \R$ defined via 
$K(\bm{x},\bm{y}):={\rm Cov}(Z(\bm{x}),Z(\bm{y}))$.
Covariance functions are symmetric and positive definite, where in this paper the term \emph{positive definite} is understood as
\begin{equation} \label{posdef}
\sum_{i=1}^n \sum_{j=1}^n {c}_i {K}(\bm{x}_i, \bm{x}_j) {c}_j \geq 0
\end{equation}
for all ${c}_i \in \mathbb{R}$, all $n \in \mathbb{N}$ and all $\bm{x}_i \in \mathbb{R}^d$.
If the inequality above is strict, 
then ${K}$ will be called strictly positive definite. 

Each symmetric 
positive definite function ${K}: \R^d \times \R^d \to \R$
defines {\em translate} functions $K(\bx,\cdot)$ on $\R^d$, for all $\bx\in\R^d$. In addition, one can define an inner product
on two translates by
\begin{equation}\label{eqHKKK}
\langle  K(\bx,\cdot), K(\by,\cdot)\rangle_{\calH(K)}
:=K(\bx, \by), \; \; \bx,\,\by\in \R^d ,
\end{equation}
in terms of $K$ itself.
This extends to all linear combinations of translates and {\em generates}, by completion, a Hilbert space $\calH(K)$ of functions on $\R^d$. 
This space is called the {\em native} space for $K$.
Notice that the Hilbert space allows for continuous point evaluations $\delta_{\bm x}\;:\;f\mapsto f({\bm x})$ via a {\em reproduction formula}
\begin{equation}\label{eqrepro}
f({\bm x})=\langle f,K({\bm x},\cdot)\rangle_{{\cal H}(K)},\;{\bm x}\in \R^d,\; \; f\in {\cal H}(K)
\end{equation}
that follows from (\ref{eqHKKK}).
Then  ${\cal H}(K)$ is called a {\em reproducing kernel Hilbert space (RKHS)} with \emph{kernel} $K$. 
In particular, the translates 
$K({\bm x},\cdot)$ lie in ${\cal H}(K)$, forming its completion and
being the Riesz representers of delta functionals $\delta_{\bm x}$.
They are central to machine learning, numerical analysis and approximation theory,
since (\ref{eqHKKK}) allows inner products in the abstract space ${\cal H}(K)$
to be explicitly computable using the kernel - the so-called \emph{kernel trick}.
See Section \ref{SecNAaAT} and
\cite{wendland2004scattered} for more detail. 
For a positive definite and stationary kernel $K$, its Fourier transform $\hat K$ can be used to recast the inner product
(\ref{eqHKKK}) on the Hilbert space ${\cal H}(K)$ by
\begin{equation}\label{eqHK}
    \langle f,g\rangle_{{\cal H}(K)}=\int_{\R^d}\frac{\hat f(\bm{\omega})
\overline{\hat g(\bm{\omega})}}{\hat K(\bm{\omega})}{\rm d}\bm{\omega},
\; \; f,g\in {\cal H}(K),
\end{equation}
up to a constant factor. 
Here, $\overline{g}$ denotes the complex conjugate of a function $g$.
Note how the spectrum of $K$ penalizes the spectrum of the
functions in ${\cal H}(K)$. 
Roughly, the Hilbert space ${\cal H}(K)$
consists of functions $f$ for which $\hat f/\sqrt{\hat K}$ is square integrable over $\mathbb{R}^d$.
The subtle connections of the Hilbert space ${\cal H}(K)$ to sample paths of Gaussian processes with covariance function $K$ will come up at many places in this paper, e.g. in Sections
\ref{Sec2}, \ref{SecQua}, \ref{PTaSP}, and \ref{sec-comp-sup}.  In this sense, kernels are important links between deterministic and probabilistic models.

A strictly positive definite kernel $K$ is called {\em stationary} 
if $K(\bm{x},\bm{y}) \equiv K(\bm{x} - \bm{y})$. According to Bochner's theorem \citep{bochner1955harmonic}, $K$ is the Fourier transform of a positive and bounded measure $F$, that is 
$$  K(\bm{x} - \bm{y}) =  \int_{\R^{d}} {\rm e}^{\mathsf{i} \left( \bm{x} - \bm{y} \, , \,  \bm{\omega} \right)} F ({\rm d} \bm{\omega}), \qquad \bm{x} , \bm{y} \in \R^d. $$
Here, $(\cdot,\cdot)$ is the inner product in $\R^d$ and $\mathsf{i}$ is the unit complex number. Fourier inversion is possible when $K$ is absolutely integrable, in which case we call the Fourier transform $\widehat{K}$ its  {\em spectral density}. 
We note that $\widehat{K}$ is nonnegative and integrable.
Furthermore, most of the paper assumes stationarity and isotropy for the covariance function,  $K$, so that 
\begin{equation}
\label{covariance}
{\rm Cov}(Z(\bm{x}),Z(\bm{y}))=K(\bm{x}-\bm{y})= \sigma^2 \varphi(\|\bm{x}-\bm{y}\|), \end{equation} for $\bm{x},\bm{y} \in \R^d$ and $\|\cdot\|$ denoting the Euclidean distance. Here, we assume $\varphi$ to be continuous with $\varphi(0)=1$. Throughout, we shall equivalently call $\varphi$ a {\em function} or a {\em correlation function}, the last as a shortcut to $\varphi(\|\cdot\|)$. Hence, the parameter $\sigma^2>0$ is the variance of $Z(\bm{x})$, for all $\bx \in \mathbb{R}^d$.
Let $\Phi_d$ denote the class of such functions $\varphi$ inducing a covariance function $K$ through the identity (\ref{covariance}) $i.e.$ $\Phi_d$ is the class of continuous isotropic correlation functions defined on $\R^d$.
Such functions have a precise integral representation according to \citet{Schoenberg1938b}, 
given by
\begin{equation}
    \label{Schoenberg}
    \varphi(x) = \int_{0}^{\infty} \Omega_d(rx) F_{d}({\rm d} r), \qquad x \ge 0,
\end{equation}
with $F_d$ being a probability measure and 
\begin{equation}
\label{Omega0} \Omega_{d}(x) = \Gamma(d/2) \left ( \frac{2}{x}\right )^{d/2-1} J_{d/2-1}(x), \qquad x \ge 0,
\end{equation}
with $\Gamma(\cdot)$ the gamma function and $J_{\nu}$ the Bessel function of the first kind of order $\nu>0$ \citep[formula 10.2.2]{Olver}. 
For a member $\varphi$ of the class $\Phi_d$, we can use that its $d$-variate Fourier transform of $\varphi(\|\bm{x}-\bm{y}\|)$ is isotropic again, and therefore
 reducible to a scalar integral formula 
\begin{equation} \label{FT}
 \widehat{\varphi}(z)= \frac{z^{1-d/2}}{(2 \pi)^{d/2}} \int_{0}^{\infty} u^{d/2} J_{d/2-1}(uz)  \varphi(u) {\rm d} u, \; z \ge 0,  
\end{equation}
defining its $d$-variate {\em isotropic spectral density}, and we assume this integral to exist. If the denominator
$(2 \pi)^{d/2}$ is omitted, the same formula holds for the inverse 
radial Fourier transform.
Throughout, we write $\Phi_{\infty}$ for $\bigcap_{d \ge 1} \Phi_d$, the class of functions $\varphi$ inducing positive definite radial functions on every $d$-dimensional Euclidean space. Hence, $\varphi \in \Phi_d$ if and only if $\varphi(\|\cdot\|)$ is a correlation function in $\R^d$. 


\section{The Mat{\'e}rn Model} \label{Sec2} 


The \emph{Mat{\'e}rn model}, ${\cal M}_{\nu, \alpha}$, is defined as \citep{Stein:1999}
\begin{equation}\label{matern}
 {\cal M}_{\nu, \alpha} (x) = \frac{2^{1-\nu}}{\Gamma(\nu)} \left ( \frac{x}{\alpha} \right )^{\nu} {\cal K}_{\nu}\left ( \frac{x}{\alpha} \right ), \qquad x \ge 0,
\end{equation}
with $\alpha > 0$ the \emph{scale} parameter, $\nu >0$ the \emph{smoothness} parameters, and ${\cal K}_{\nu}$ a modified Bessel function of the second kind of order 
 $\nu$ \citep[][9.6.22]{Abra:Steg:70}. 
 It can be verified that ${\cal M}_{\nu,\alpha}(0)=1$, so that \eqref{matern} is a correlation function.
 Arguments in \citet[][p48]{Stein:1999} show that ${\cal M}_{\nu,\alpha}$ belongs to the class $\Phi_{\infty}$.
The function $\sigma^2 {\cal M}_{\nu,\alpha}$ will be termed \emph{Mat{\'e}rn covariance function}, and $\sigma^2>0$ will denote the variance of the associated Gaussian random field.

The importance of the Mat{\'e}rn class stems from the parameter $\nu$ that controls the differentiability  of
the sample paths
 of the associated Gaussian field. Specifically, for
any positive integer $k$,
 the
sample paths 
of 
 a Gaussian field $Z$ on $\R^d$
 with Mat{\'e}rn correlation function are  $k$-times mean square differentiable (in any direction) if and only if $\nu > k$. Also,
 a rescaled version of the Mat{\'e}rn correlation function converges to the Gaussian or squared exponential kernel
as $\nu \rightarrow \infty$, that is 
\begin{equation}\label{eqGauLim}
{\cal M}_{\nu, \alpha/(2 \sqrt{\nu})}(x) \xrightarrow[\nu \to \infty]{}  \exp(- x^2/\alpha^2), \qquad x \ge 0, 
\end{equation}
with convergence being uniform on any compact set of $\R^d$.
For this reason, the parametrisation ${\cal M}_{\nu, \alpha/(2 \sqrt{\nu})}$ is sometimes also adopted \citep{williams2006gaussian}. 

When $\nu=k+ 1/2$, for $k$ a nonnegative integer,  the Mat{\'e}rn correlation function simplifies into the product of a negative exponential correlation function with a polynomial of order $k$. 
For instance, ${\cal M}_{1/2,1}(x)= \exp(-x)$ and ${\cal M}_{3/2,1}(x)= \exp(-x) (1+x)$.
In general,  
\begin{equation}
{\cal M}_{k+1/2,1}(x)= \exp(-x) \sum_{i=0}^k \frac{(k+i)!}{2k!}  \binom{k}{i}    (2x)^{k-i}
\end{equation}
for $k \in \mathbb{N}_0$.
This simple algebraic form for the \M correlation functions has undoubtedly contributed to the widespread popularity of the \M model.

Now we are in a position to explore in detail the many faces of the \M model.
Section \ref{Sec3} discusses maximum likelihood estimation, Kriging prediction, and the screening effect,
while Section \ref{Sec4} explores an SPDE characterisation of the \M model. 
Section \ref{secModeling} discusses the Mat{\'e}rn model as a building block to more sophisticated models, while Section \ref{sec5} views the scientific landscape through the lens of the \M model, with special emphasis on numerical analysis, probability theory and machine learning. 
Section \ref{sec6} introduces some recently developed alternatives and generalisations of the Mat{\'e}rn model, while
Section \ref{sec7} compares these alternative models in terms of estimation, prediction, and the screening effect.

\section{Estimation and Prediction with the Mat{\'e}rn Model} \label{Sec3}




Let $D\subset \R^d$ be  a subset of $ \R^d$. 
Consider a set $X_n = \{\bm{x}_1,\dots,\bm{x}_n\}$  of (distinct) locations in $D$, at which values $\boldsymbol{Z}_n=(Z(\boldsymbol{x}_1),\ldots,Z(\boldsymbol{x}_n))^{\top}$
of the Gaussian random field $Z$, defined in Section \ref{background}, are observed. 
An important problem concerns the \emph{prediction} of values $Z(\bm{x}_0)$ at an unobserved location $\boldsymbol{x}_0\in D \setminus X_n$.
Then an especially natural predictor for $Z(\bm{x}_0)$ is 
\begin{equation}\label{blup}
\widehat{Z}_{n} =\boldsymbol{c}_n^\top 
 \bR_{n}^{-1}\boldsymbol{Z}_n
\end{equation}
with the vector 
$[\boldsymbol{c}_n]_i= K(\bm{x}_0,\bm{x}_i)$ and 
the {\em kernel matrix} 
$[\bR_{n}]_{i,j} = K(\bm{x}_i,\bm{x}_j)$. 
The predictor \eqref{blup} can be motivated from multiple directions.
Classically, \eqref{blup} is motivated as the 
best linear unbiased predictor (BLUP) for $Z(\bm{x}_0)$, and is often referred to as the {\em simple Kriging} 
predictor of $Z(\bm{x}_0)$ \citep{cressie1990origins}.
From a modern perspective, where the role of unbiased estimation is increasingly questioned, 
we can motivate this choice using alternative optimality properties, including:
\begin{enumerate}
\item it is the 
expectation 
of $Z(\bm{x}_0)$ conditionally on the realisation $\bm{Z}_n$;
\item it is the optimal estimate (i.e. the Bayes act) for $Z(\bm{x}_0)$ based on the data-set $\bm{Z}_n$, under squared error loss \citep[][Section 13.3]{novak2008tractability};
\item it yields the minimal RKHS norm interpolant of the data evaluated at $\bm{x}_0$, by Section \ref{SecNAaAT};
\item it is the algorithm for approximating $Z(\bm{x}_0)$ from $\bm{Z}_n$ that minimises the worst case error in the sense of information-based complexity \citep[][Section 10.2]{novak2008tractability} 
and approximation theory (see Section \ref{SecNAaAT}),
\end{enumerate}
to name but a few. 
The Mat\'{e}rn model provides a natural setting to study the performance of \eqref{blup} if we suppose $Z$ to have a stationary isotropic covariance function $\sigma^2 {\cal M}_{\nu,\alpha}$. 
The crucial question of how to select suitable values for the parameters $\sigma$, $\alpha$, $\nu$ will be considered first, in Section \ref{seclik}, and then the performance of \eqref{blup} will be studied in Section \ref{pprr}.
The possibility of a direct extension of the \M model to more general domains, such as manifolds and graphs, is discussed in Section \ref{subsec: manifolds and graphs}.

\subsection{Estimation Using Maximum Likelihood} \label{seclik}

Maximum likelihood 
(ML) and similar estimation methods are popular in this setting
due to the availability of practical (inc. gradient-based) numerical methods for computation and the classical theory that underpins ML. 
On the other hand, implicit in the use of ML is that the statistical model is well-specified, and this judgement must be made on a case-by-case basis.
To limit scope, we focus on ML estimation in the sequel. 
Our aim is to understand when the parameters of the Mat\'{e}rn model can be consistently estimated from data, and to understand  the asymptotic distribution of  the  ML estimator.
To this end, recall that the Gaussian log-likelihood function is 
\begin{eqnarray}\label{eq:17}
\hspace{20pt} \mathcal{L}_{n}(\boldsymbol{\theta}) = -\frac{1}{2} \left( \log(|\sigma^2 \bR_{n})|) + \frac{1}{\sigma^{2}}\boldsymbol{Z}_n^{\top}\bR_{n}^{-1}\boldsymbol{Z}_n \right),
\end{eqnarray}
up to an additive constant, with $\boldsymbol{\theta}=(\nu,\alpha, \sigma^{2})$.
 The ML estimator is defined as 
\begin{equation}
\widehat{\boldsymbol{\theta}}_{n}=\underset{\boldsymbol{\theta} \in \R^3_{+} }{\textrm{argmax}}\,\mathcal{L}_{n}(\boldsymbol{\theta}) . \label{eq: max like}
\end{equation}
The  ML estimate for the variance parameter can be  computed in closed-form as $\widehat{\sigma}^2_n=\boldsymbol{Z}_n^{\top} \bR^{-1}_n \boldsymbol{Z}_n/n$; plugging this expression into \eqref{eq:17} reduces the numerical problem to optimisation of a so-called \emph{concentrated likelihood} over $\mathbb{R}_+^2$.
 However, maximizing the log-(concentrated) 
 likelihood 
 requires a  nonlinear optimisation problem to be solved, for which numerical methods must be used; see Section \ref{subsec: approx like}.

The performance of ML estimation has been studied principally in two different asymptotic limits.
Under \emph{fixed domain asymptotics}, 
 the sampling domain $D$ is bounded and the set of sampled locations
$X_n$ becomes increasingly dense in $D$. 
Under \emph{increasing domain asymptotics}, the domain $D$ grows with the number $n$ of observed data, and the distance between any two sampled locations is bounded away from zero.
\citet{Zhang:Zimmerman:2005} note that the peformance of the ML estimator can be quite different under these two frameworks, as will now be discussed. 

\subsubsection{Increasing Domain Asymptotics.}\label{ida}

\citet{Mardia:Marshall:1984} make use of increasing domain asymptotics 
to establish, under mild regularity conditions, that the ML estimator is \emph{strongly consistent}, meaning that
 $\widehat{\boldsymbol{\theta}}_{n}\stackrel{a.s.}{\longrightarrow} \boldsymbol{\theta}_0$ for the {\em true} parameter $\boldsymbol{\psi}_0$. 
Furthermore, they establish that the ML estimator is \emph{asymptotically normal}, meaning that
 \begin{equation}
 \bm{F}^{1/2}(\boldsymbol{\theta}_0)(\widehat{\boldsymbol{\theta}}_{n}-\boldsymbol{\theta}_0) \stackrel{d}{\longrightarrow} \mathcal{N}(\boldsymbol{0},\bm{I})
 \end{equation}
  where
  $\bm{F}(\boldsymbol{\theta})=-E[\mathcal{L}^{''}_{n}(\boldsymbol{\theta})]$ is the Fisher information matrix, whose entries are
$$
F(\boldsymbol{\theta})_{i,j}=
\frac {1}{2}{\rm tr} \left(\frac{\mathrm{d} \bm{\Sigma}_n}{\mathrm{d} \boldsymbol{\theta}_i} \bm{\Sigma}_n^{-1}\frac{\mathrm{d} \bm{\Sigma}_n}{\mathrm{d} \boldsymbol{\theta}_j} \bm{\Sigma}_n^{-1}\right) ,
$$
and $\bm{\Sigma}_{n}=\sigma^2 \bR_n$.
Although our focus is on the \M model, we note that these kind of asymptotic results hold for any parametric correlation function obeying particular regularity   conditions that   are stated in terms of  eigenvalue conditions on the correlation matrix and its derivatives \citep{Mardia:Marshall:1984}, thought these may not be easy to verify in general (see for instance \citet{Shaby:Ruppert:2012}, for the exponential case).
Generally speaking,  as long as the 
spatial extent of the sampling region is large compared with
the range of dependence of the random field, increasing-domain asymptotics provide a very
accurate description of the behavior of the ML estimate 
\citep{Zhang:Zimmerman:2005,Shaby:Ruppert:2012,Shaby:Kaufmann:2013}.

\subsubsection{Fixed Domain Asymptotics.}
 
\citet{Zhang:2004} considered ML estimation for the Mat{\'e}rn model under fixed domain asymptotics, proving that when the smoothness
parameter $\nu$ is known and fixed, none of the parameters $\sigma^2$ and $\alpha$ can be estimated
consistently when $d=1, 2, 3$. Instead, only the parameter 
\begin{equation} \label{micro_matern}
\text{micro}_{{\cal M}}= {\sigma^2}/{ \alpha^{2\nu}},
\end{equation}
 sometimes called \emph{microergodic} parameter \citep{Zhang:Zimmerman:2005,Stein:1999}, can be
consistently estimated.
This is a consequence of the equivalence of the two corresponding Gaussian measures,
that we denote with $P(\sigma^2_i{\cal M}_{\nu,\alpha_i})$, with $i=0,1$.
In particular, for any
bounded infinite set $D\subset \R^d$, $d=1, 2, 3$,
$P(\sigma^2_0{\cal M}_{\nu,\alpha_0})$ is equivalent to  $P( \sigma^2_1{\cal M}_{\nu,\alpha_1})$ on the paths of $Z(\boldsymbol{x}), \boldsymbol{x} \in D$, if and only if 
\begin{equation}
 \label{matern_equivalence}
 \sigma_0^2 / \alpha_0^{2\nu}=\sigma_1^2/ \alpha_1^{2\nu}.
\end{equation}
In contrast, for $d \geq 5$, \citet{anderes2010} proved  the  orthogonality of two Gaussian measures
with different Mat{\'e}rn covariance functions
and hence, in this case, all the parameters can be consistently estimated under fixed-domain asymptotics. The case $d = 4$ has been recently studied in \citet{bolin3}.

Asymptotic results associated with  ML estimation of the microergodic parameter, again for a fixed known smoothness parameter $\nu$, can be found in \citet{Zhang:2004}, and later on in
\citet{Shaby:Kaufmann:2013}. 
In particular, 
for a zero mean Gaussian field  defined on a  bounded infinite set $D\subset \R^d$, $d=1, 2, 3$, 
with a Mat{\'e}rn  covariance function $\sigma^2_0{\cal M}_{\nu,\alpha_0}$
the ML estimator $\hat{\sigma}_{n}^{2}/\hat{\alpha}_{n}^{2\nu}$ of the microergodic parameter 
is strongly  consistent,  $i.e.$,  $$\hat{\sigma}_{n}^{2}/\hat{\alpha}_{n}^{2\nu}\stackrel{a.s.}{\longrightarrow} \sigma_{0}^{2}/\alpha_{0}^{2\nu},$$ and its  asymptotic distribution is given by 
$$\sqrt{n}(\hat{\sigma}_{n}^{2}/\hat{\alpha}_{n}^{2\nu}-\sigma_{0}^{2}/\alpha_{0}^{2\nu})\stackrel{d}{\longrightarrow} \mathcal{N}(0,2(\sigma_{0}^{2}/\alpha_{0}^{2\nu})^2).$$
Generally speaking,
when the  range of dependence of the random field is large with respect to the spatial extent of the sampling region,  fixed domain  asymptotics provide a very
accurate description of the behavior of the ML estimate of the microergodic parameter
\citep{Shaby:Kaufmann:2013}.
Extensions of these results to the case where $Z$ is observed with Gaussian errors can be found in \citet{tang_ban2021}, while results for a space-time version of the  Mat{\'e}rn model  
can be found in \citet{LG} and \citet{faouzi2022space}.
Finally we highlight that the efficient estimation of the microergodic parameter assuming the smoothness parameter unknown is still an open problem; some promising results in this direction can be found in \citet{AOS2077}.

\subsection{Prediction and the Screening Effect}
\label{pprr}

The equivalence of Gaussian measures within the 
Mat{\'e}rn class has consequences for prediction of $Z(\bm{x}_0)$ at an unobserved location $\bm{x}_0 \in D \setminus X_n$; these consequences will now be discussed.
In what follows, $\nu$ is supposed known and fixed, and we consider the setting where $\sigma$ and $\alpha$ are \emph{misspecified}.
That is, we suppose $Z$ is a Gaussian field with Mat{\'e}rn covariance $\sigma_0^2 {\cal M}_{\nu,\alpha_0}$, and we consider the performance of the predictor \eqref{blup} when a Mat{\'e}rn model $\sigma_1^2 {\cal M}_{\nu,\alpha_1}$ is used. 
This situation is typical, since the true parameters $\sigma_0$ and $\alpha_0$ of the data-generating process will be unknown in general.
Our theoretical setting will be fixed domain asymptotics.

Note, first, that \eqref{blup} does not depend on the value of $\sigma_1$, but does depend on the value of the parameter $\alpha_1$ (and the parameter $\nu$, but this parameter is fixed).
This dependence will be emphasised using the notation $\boldsymbol{c}_n(\alpha_1)$ and  $\bR_{n}(\alpha_1)$.
Under the Gaussian measure $P(\sigma^2_0{\cal M}_{\nu,\alpha_0})$ associated with the {\em true} model $\sigma^2_0{\cal M}_{\nu,\alpha_0}$, the mean squared error of the predictor $\widehat{Z}_{n}(\alpha_1)$ is given by
\begin{align*}
&\text{{\sc var}}_{\alpha_0,\sigma^2_0}\left[\widehat{Z}_{n}(\alpha_1)-Z(\boldsymbol{x}_0)\right] \\ 
&=\sigma_0^2\Big(1-2\boldsymbol{c}_n(\alpha_1)^{\top
} \bR_{n}(\alpha_1)^{-1}
\boldsymbol{c}_n(\alpha_0) \\ 
&\hspace{30pt} + \boldsymbol{c}_n(\alpha_1)^{\top} \bR_{n}(\alpha_1)^{-1} \bR_{n}(\alpha_0) \bR_{n}(\alpha_1)^{-1}\boldsymbol{c}_n(\alpha_1)\Big) ,
\end{align*}
while if there is no misspecification then the previous expression reduces to
\begin{align}\label{msetrue}
& \text{{\sc var}}_{\alpha_0,\sigma_0^2}\big[\widehat{Z}_{n}(\alpha_0)-Z(\boldsymbol{x}_0)\big]\\
 \nonumber   & \hspace{40pt} = \sigma_0^2\big(1-\boldsymbol{c}_n(\alpha_0)^\top\bR_{n}^{-1}(\alpha_0)\boldsymbol{c}_n(\alpha_0)\big).
\end{align}
Under regularity conditions, and for fixed domain asymptotics, \citet{Stein:1988} shows that both  asymptotically efficient prediction and asymptotically correct estimation of prediction variance
 hold when the two  Gaussian measures  $P(\sigma^2_i{\cal M}_{\nu,\alpha_i})$, $i=0,1$ are equivalent, {\em i.e.} 
 \eqref{matern_equivalence}.
 Specifically, 
 \begin{equation}\label{kauf3_1} \frac{\text{{\sc var}}_{\sigma^2_0,\alpha_0}\bigl[\widehat{Z}_{n}(\alpha_1)-Z(\boldsymbol{x}_0)\bigr]}{\text{{\sc var}}_{\sigma^2_0,\alpha_0}\bigl[\widehat{Z}_{n}(\alpha_0)-Z(\boldsymbol{x}_0)\bigr]}\stackrel{a.s.}{\longrightarrow}1
        \end{equation}
and 
 \begin{equation}\label{kauf3_3} \frac{\text{{\sc var}}_{\sigma^2_1,\alpha_1}\bigl[\widehat{Z}_{n}(\alpha_1)-Z(\boldsymbol{x}_0)\bigr]}{\text{{\sc var}}_{\sigma^2_0,\alpha_0}\bigl[\widehat{Z}_{n}(\alpha_1)-Z(\boldsymbol{x}_0)\bigr]}\stackrel{a.s.}{\longrightarrow}1.
      \end{equation}
The implication of  (\ref{kauf3_1}) is that, under the true model,  if the correct value of $\nu$ is used,
any value of $\alpha_1$  will give asymptotic efficiency. The  implication of  (\ref{kauf3_3})
is stronger and guarantees  that using the misspecified predictor 
under the correct and misspecified models  is asymptotically equivalent from mean squared error point of view.
Note that these kind of results does not consider the uncertainty associated with the covariance parameters
of the misspecified model.
\citet{Shaby:Kaufmann:2013} show that  (\ref{kauf3_3}) still holds by considering the ML estimator of the variance $\hat{\sigma}^2_n=\boldsymbol{Z}_n^\top \bR^{-1}_n(\alpha_1) \boldsymbol{Z}_n/n$ in place $\sigma^2_1$.

Conditions of equivalence of two Gaussian measures based on a space-time \cite{LG} and bivariate \cite{frbev} version of the Mat{\'e}rn  model have also been established.
Next, we consider a practically important aspect of prediction; the co-called screening effect.

\vspace{0.3cm} 
 {\em Screening Effect.}  The \emph{screening effect} refers to the phenomenon where the predictor \eqref{blup} depends almost exclusively on those observations that are located nearest to the predictand \citep{stein1}. 
 As such, the screening effect is an important tool that can be used to mitigate the computational burden of evaluating \eqref{blup} in the presence of big datasets. 
 This issue has traditionally been an important subject in geostatistics \citep{matheron63, matheron65, matheron71, chiles}. 
Indeed, \citet{matheron63, matheron65}, in the School of Geostatistics at the Ecole des Mines, developed a first formalisation of screening effect, referring to situations where the observations located far from the predictand receive a zero kriging weight. Matheron's definition has a direct connection with the Markov property on the real line, which happens when kriging is performed under the exponential model (indeed, ${\cal M}_{1/2,\alpha}$). 


M. Stein 
\citep{Stein:1999, stein1, stein2, stein2015does} 
adopts an alternative definition of the screening effect that will now be described. 
Let 
$Z$ be a mean-square continuous, zero mean and weakly stationary Gaussian random field on $\R^d$.
 Let $e(X_n)$ be the error of the predictor \eqref{blup} of $Z(\bm{x}_0)$ based on $\bm{Z}_n$. 
 Two choices for the set $X_n$ of observation locations will be considered, and to this end we let $F_{\epsilon}, N_{\epsilon}$ be sets, indexed by $\epsilon>0$, such that $N_{\epsilon}$ contains the nearest observations to the predictand, and $F_{\epsilon}$ the furthest observations. 
 Then \citet{stein1} says that $N_{\epsilon}$ {\it asymptotically screens out} $F_{\epsilon}$ when
\def\e1{\epsilon}
\begin{equation}
\label{def-screening} \lim_{\e1 \downarrow 0} \frac{\mathbb{E} \; e (N_{\e1} \cup F_{\e1})^2}{\mathbb{E} \; e (N_{\e1})^2} = 1.
\end{equation}
A thorough discussion of the implications of this definition can be found in \citet{porcu2020stein}, where nontrivial differences between fixed domain and increasing domain asymptotics are reported. 



The spatial configuration of the sampling point $X_n$ determines whether the screening effect will hold. 
\citet{porcu2020stein} refer to a {\em regular scheme} as one for which $F_{\e1} = \{ \e1 (\bm{x}_0+j) \}$, for $j \in \mathbb{Z}^d$ and $N_{\e1}$ being the restriction of $F_{\e1}$ to some fixed region with $\bm{x}_0$ in its interior, assuming $\bm{x}_0 \notin \Z^d $.
For regular schemes, \citet{stein1} established (\ref{def-screening}) whenever {\em the spectrum $\widehat{K}$ varies regularly at infinity \citep{bingham} in every direction with a common index of variation} \citep[quoted from][]{porcu2020stein}.  
However, this condition may not be useful for space-time processes, where differentiability properties in the space and time coordinates are not necessarily identical. 
To overcome such a problem, we instead consider an {\em irregular scheme}: for $\bm{x}_1, \ldots, \bm{x}_n$ being distinct nonzero elements of $\R^d$, $\bm{y}_1, \ldots, \bm{y}_N$  distinct elements of $\R^d$, $\bm{x}_0=\bm{0} \in \R^d$ and $\bm{y}_0 \in \R^d$ being nonzero, we have $N_{\e1} =\{\e1 \bm{x}_1,\ldots,\e1 \bm{x}_n \}$ and $F_{\e1} =\{\bm{y}_0 +\e1\bm{y}_1,\ldots, \bm{y}_0 +\e1\bm{y}_N\}$. 
The \emph{Stein hypothesis} \citep[termed in][]{porcu2020stein}
\begin{equation} \label{condition1}
\forall R > 0, \quad \lim_{\|\bm{\omega}\| \to \infty } \sup_{\| \bm{\tau}\|<R} \bigg | \frac{\widehat{K}(\bm{\omega}+\bm{\tau} )}{\widehat{K}(\bm{\omega})}  -1 \bigg |=0,
\end{equation} 
provides a sufficient condition for the screening effect in this setting (under some mild additional conditions on $\widehat{K}$ and $N_\e1$),
which can be verified in dimensions $d=1$ and $d=2$ for mean-square continuous but non-differentiable random fields, for some specific designs $N_\e1$ \citep{stein2}. 
The \M model with $K = {\cal M}_{\alpha,\nu}$ admits
 a simple expression for its spectrum \citep[][11.4.44]{Abra:Steg:70}:  
\begin{equation} \label{stein1}
\widehat{{\cal M}}_{\nu,\alpha}(z)= \frac{\Gamma(\nu+d/2)}{\pi^{d/2} \Gamma(\nu)}
\frac{ \alpha^d}{(1+\alpha^2z^2)^{\nu+d/2}}
, \; \; z \ge 0 ,
\end{equation}
from which \eqref{condition1} can be 
verified.

The screening effect can thus be established for the Mat{\'e}rn model, under both regular and irregular schemes, justifying the use of ``local'' approximations to the predictor \eqref{blup}.


\subsection{Mat{\'e}rn on Manifolds and Graphs}
\label{subsec: manifolds and graphs}


Let $M$ be a general manifold.
A pragmatic question is whether the Mat{\'e}rn correlation function \eqref{matern} can be composed with a suitable metric $g$, defined on the manifold, to preserve positive definiteness over $M$. 
For the case of the sphere, a natural metric is the geodesic distance; the length of the arc connecting any pair of points located over the spherical shell. 
For this metric, $(x,y) \mapsto {\cal M}_{\nu,\alpha}(g(x,y))$ is a correlation function only for $0 < \nu \le 1/2$ \citep{gneiting2013}. 
This limitation is emphasised in \citet{alegria2021f}, who propose the ${\cal F}$ family, a model that is valid on the sphere, and having the same properties as the Mat{\'e}rn function in terms of mean-square differentiability and fractal dimension.
The Mat{\'e}rn function on other general manifolds has been studied by \citet{li2021inference}. \citet{guinness2016isotropic} propose a spectral expansion to define a covariance function that mimics the Mat{\'e}rn model, but this construction is criticised by \citet{lindgren2022spde} as being incorrect as the spectral expansion does not reproduce the same properties of the \M model.

Unfortunately, it seems that the limited applicability of the Mat{\'e}rn model on any space that is not a flat surface extends to more abstract settings as well. An elegant isometric embedding argument in \citet{anderes2020isotropic} proves that the restriction $0 < \nu \le 1/2$ is required when the input space is a graph with Euclidean edges.
A more general argument in \citet{menegatto2020gneiting} proves that the same restriction is inherited for a general quasi metric space endowed with a geodesic metric. The notable effort by \citet{bolin2020rational} provides a model that is once differentiable over metric graphs. 
It is reasonable to conclude that some form of the SPDE approach, which we discuss next in Section \ref{SecSPDE}, is needed in general to extend the \M model to a general manifold.



\section{The Mat{\'e}rn Model in Computational Statistics } 
\label{Sec4}

This section explores the interaction of the \M model with computational statistics, starting with numerical methods for \emph{implementation} of the \M model (Sections \ref{SecSDE}, \ref{SecSPDE} and \ref{subsec: approx like}), and then turning to uses of the \M model to \emph{facilitate} numerical computation itself (Section \ref{SecQua}).


\subsection{Implementation as an SDE}
\label{SecSDE}

The \M model admits a \emph{state space} representation as an SDE, which enables efficient computational techniques from the signal processing literature to be employed for simulation, estimation and prediction.
Indeed, focusing on dimension $d=1$, and letting 
$$
\mathbf{Z}(x) = (Z,\mathrm{d}Z/\mathrm{d}x,\dots,\mathrm{d}^k Z/\mathrm{d}x^k) ,
$$
the \M model $\mathcal{M}_{\nu,\alpha}$ with $\nu = k + 1/2$ admits the characterisation
\begin{align*}
    \mathrm{d}\mathbf{Z} = \left( \begin{array}{cccc} 0 & 1 & & \\ & \vdots & \vdots & \\ & & 0 & 1 \\ - a_0 & - a_1 & \dots & - a_{k-1}  \end{array} \right) \mathbf{Z} \; \mathrm{d}x + \left( \begin{array}{c} 0 \\ \vdots \\ 0 \\ 1 \end{array} \right) \mathrm{d}\mathcal{W}
\end{align*}
where $a_i = {}_{k+1}\text{C}_i \cdot \alpha^{-k-1+i}$, the ${}_\cdot\text{C}_\cdot$ are binomial coefficients, and $\mathcal{W}(x)$ represents a zero-mean white noise process on $x \in \mathbb{R}$ \cite{hartikainen2010kalman}.
The advantage of state space formulations is that both estimation and prediction can be performed in a \emph{single pass} through the data, at linear $O(n)$ cost, using familiar Kalman updating equations as described in \citet{sarkka2013SDE} and in further detail in Chapter 6 of \citet{hennig2022probabilistic}.
Similar characterisations for higher dimensions, including spatio-temporal versions of the \M model, can be found in \citet{sarkka2013SDE}, though we note these retain linear complexity only in the number of time steps; complexity is cubic in the size of the spatial grid.
The SPDE approach can offer a solution in this respect, and we discuss this next.

\subsection{Implementation as an SPDE}
\label{SecSPDE}



A major reason for the continued popularity of the \M model is the availability of efficient and scalable numerical methods for simulation, due in large part to \citet{Lindgren}.
These authors consider the SPDE 
\begin{equation}
\label{SPDE} 
(\alpha^{-2} -\Delta  )^{\gamma/2} Z(\bm{x}) = {\cal W}(\bm{x}), \qquad \bm{x} \in \mathbb{R}^d,
\end{equation}
where $\alpha>0$, $\Delta$ is the Laplacian, 
and ${\cal W}$ is a Gaussian white noise on $\R^d$, so that $\text{Cov}\left ( {\cal W}(A_1),{\cal W}(A_2) \right )= |A_1 \cap A_2|$, where $A_i$ are subsets of $\R^d$, $i=1,2$, and where $|\cdot|$ is the volume integral. 
\citet{whittle} and \citet{whittle1963stochastic} proved that the solution to \eqref{SPDE} is a Gaussian field with Mat{\'e}rn covariance $\sigma^2 \mathcal{M}_{\nu,\alpha}$ with parameters $\alpha$ (as before) and
$$
\sigma^2 = \frac{\Gamma(\nu) \alpha^{2 \nu}}{\Gamma(\nu + d/2) (4 \pi)^{d/2}} , \qquad \nu=\gamma-d/2 .
$$
This perspective offers two insights; first, tools developed for the numerical approximation of SPDEs can be brought to bear on the \M model, and second, there is a clear path to generalise the definition of the \M model to 
any (planar or non planar) manifold on which the analogous SPDE may be defined. 
(For example, \citet{jansson2022surface} take this perspective to generalise the \M model to the sphere $\mathbb{S}^d$.)

To provide a computationally convenient approximation to (\ref{SPDE}), \citet{Lindgren} 
considered the weak solution to (\ref{SPDE}) 
and approximation of the weak solution using basis functions with compact support over a compact domain $\Omega \subset \mathbb{R}^d$ (specifically, a \emph{Galerkin} approximation using finite element basis functions was used). 
As a result, the authors establish a formal route to approximation of the random field $Z$ with a \emph{Gauss--Markov} random field having a \emph{sparse} precision matrix. 
Sparse matrix algebra enables fast simulation of realisations from the \M random field, and fast evaluation of the likelihood \eqref{eq:17} (albeit not fast evaluation of the gradient of the likelihood).

The choice of domain $\Omega$ introduces boundary effects which must be carefully mitigated.
\citet{khristenko2019analysis,brown2020semivariogram} provide a solution for the case where $\gamma$ is an integer; the non-integer case is considered in \citet{bolin2020rational}. 
The extension of the Mat{\'e}rn field based on SPDEs to space-time is provided by \citet{cameletti2013spatio} and subsequently by \citet{bakka2020diffusion,clarotto2022spde}, while the multivariate Mat{\'e}rn case has been explored in \citet{bolin2016multivariate}. 
Alternative approximations based on Galerkin methods on manifolds have been provided by \citet{lang2021galerkin}. 
An interesting approach that allows working on manifolds with huge datasets is proposed by \citet{pereira2022geostatistics}. 
 The interest in this literature is dual. 
 On the one hand, the technical aspects related to the finite dimensional representation of Gaussian random fields are extremely interesting {\em per se}. 
 On the other hand, this group of authors is actually driven by providing tools for efficient computation. This is witnessed by the relevant existing packages (R-INLA, inlabru, and rSPDE for instance) and we refer to the review of \citet{lindgren2022spde}.

\citet{sanz2022finite} attempt to explain the trade-off between accuracy and scalability in numerical approxmation of the \M model. 
Recall that, in the SPDE approach \cite{Lindgren}, $Z$ in (\ref{SPDE}) is numerically approximated using a Gaussian process 
\begin{equation}
    \label{delta}
    Z_{\delta}(\bm{x})= \sum_{k=1}^{n_{\delta}} \omega_k \epsilon_k(\bm{x}), \qquad \bm{x} \in \Omega,
\end{equation}
where $\epsilon_k$ are finite element basis functions and the vector $\bm{\omega}=(\omega_1, \ldots, \omega_{n_{\delta}})^{\top}$ is multivariate Gaussian with zero mean and with a sparse precision matrix. 
The accuracy of the approximation $Z_\delta$ is dependent on (a) the compact support of the finite elements basis functions, (b) boundary effects due to the domain $\Omega$,
and (c) by the mesh width $\delta$ that determines the cardinality $n_{\delta}$ in (\ref{delta}). 
Most of the earlier literature has considered (\ref{delta}) with $n_{\delta}$ proportional to the sample size $n$ of the dataset being modelled.  
\citet{sanz2022finite} adopt a fixed domain asymptotic approach to explain when $n_{\delta} \ll n$ might be a legitimate strategy. 
To do so, they consider Gaussian process regression and work under the framework of Bayesian contraction rates. 
Their results provide justification for specific scalings of $n_\delta$ with $n_\delta = o(n)$, provided that the smoothness $\nu$ is sufficiently high.


A different path to SPDE and Gauss--Markov random fields was recently taken in \citet{sanz2022spde}, who adopt graph-based discretisations of SPDEs.
This approach can be well-suited to working with discrete and unstructured point clouds, such as in machine learning tasks where the data belong to an implicitly defined low-dimensional manifold. 
A second advantage of this approach is that an explicit triangulation of the domain is not required.

\subsection{Approximate Likelihood and the Mat{\'e}rn Model}
\label{subsec: approx like}


In estimating the parameters of the \M model using ML \eqref{eq: max like}, numerical optimisation is required.
Although generic optimisation routines can be used, an often better approach is to first construct an accurate-but-cheap approximation to the likelihood, which can then be more readily maximised. 
Indeed, approximate likelihoods are essential when dealing with large datasets, since the evaluation of (\ref{eq:17})  
requires computing the inverse and the determinant of the correlation matrix, usually via the Cholesky decomposition at complexity $O(n^3)$ and storage cost $O(n^2)$. 

Perhaps the most successful approximation is \emph{Vecchia’s method} \citep{vecchia1988}, which has attracted a remarkable amount of attention in recent times \citep[inc.][]{Stein:2005ba,datta2016hierarchical,datta2016nonseparable,guiness2018,datta2021}. 
The Vecchia approximation can be used  with any correlation model and its basic idea is  is to replace (\ref{eq:17}) with a product of Gaussian conditional distributions, in which each conditional distribution involves only a small subset of the data. 
This approximation requires that the data are \emph{ordered} and the number $m$ of `previous' data on which to condition is to be specified. 
Generally, larger $m$ entails more accurate and computationally expensive approximation, 
while the choice of ordering  affects the accuracy of the approximation \citep{guiness2018}.
The Vecchia method provides a sparse 
approximation to the Cholesky factor of the precision matrix, such that the approximate likelihood can be computed in $O(nm^3)$ time and with $O(nm^2)$ storage cost. 
See the recent review of \citet{ka2021} for further detail.
The Vecchia likelihood  can be viewed as a specific instance of a more general class of estimation methods called quasi- or composite likelihood \citep{Lindsay:1988,varin2011overview} that have been widely used for the estimation of Gaussian fields with the 
Mat{\'e}rn model \citep{Eidsvik:Shaby:Reich:Wheeler:Niemi:2013,bevilacqua2015comparing,BACHOC2019}.


An alternative method of mitigating the computational burden  of ML estimation is \emph{covariance tapering} \citep{furrer2006covariance}.
The basic idea is to multiply  the Mat{\'e}rn model  with a
compactly supported correlation function, resulting in a `modified' Mat{\'e}rn model with compact support. 
This induces sparseness in the associated covariance matrix, so that algorithms for sparse matrices can be exploited
for a computationally efficient evaluation of the Cholesky decomposition \citep{furrer2006covariance}.
However, some authors \citep{BFFP,bevilacqua2022unifying} suggest that  tapering might be an obsolete approach in view of the fact that flexible compactly supported models that include the Mat{\'e}rn model as a special case have been recently proposed; see Section \ref{sec7}.
A comprehensive review of the likelihood approximations is beyond the scopes of this paper, so we refer the reader to 
\citet{sun2012geostatistics} and \citet{popurri}
for further detail.

\subsection{The Mat\'{e}rn Model \emph{for} Bayesian Computation}\label{SecQua}

In the last decade there has been increasing interest in the use of kernel methods for solving PDEs.
Consider a system
\begin{align*}
\mathcal{A} u & = f \qquad \text{in } \Omega \\
\mathcal{B} u & = g \qquad \text{on } \partial\Omega
\end{align*}
specified by a differential equation involving $\mathcal{A}$ and $f$, and initial or boundary conditions specified by $\mathcal{B}$ and $g$.
Dating back at least to \citet{fasshauer1996solving} in the deterministic setting, and reinterpreted through a Bayesian lens by authors such as \citet{cockayne2019bayesian}, one can seek an approximation to the strong solution $u : \Omega \rightarrow \mathbb{R}$ by modelling $u$ as \emph{a priori} a Gaussian random field and conditioning that field to satisfy the differential equation at locations $\{\bm{x}_1,\dots,\bm{x}_m\} \subset \Omega$ and satisfy the boundary conditions at locations $\{\bm{x}_{m+1},\dots,\bm{x}_n\} \subset \partial \Omega$. 
The conditional mean of this process coincides with the \emph{symmetric collocation} method introduced by  \citet{fasshauer1996solving}, which we return to in Section 
\ref{SecNAaAT}, while the conditional variance provides probabilistic uncertainty quantification for the solution, expressing the uncertainty that remains as a result of using only a finite computational budget.
To implement these methods, one requires a Gaussian process whose sample paths possess sufficient regularity for the operation of conditioning on the derivative $\mathcal{A}u$ to be well-defined.
On the other hand, assuming excessive smoothness could lead to over-confident uncertainty quantification.
One therefore requires a kernel with customisable smoothness, which can be adapted to the differential equation at hand.
The Mat\'{e}rn class satisfies this requirement, but is not alone in doing so; we continue discussion of this point in Section \ref{surr}.

A specific PDE that has received considerable recent attention in the Bayesian statistical community is the \emph{Stein equation}, for which $\mathcal{A}u = c + p^{-1} \nabla \cdot (p \nabla u)$, where $p$ is the probability density function of a posterior distribution of interest, $f$ is a function whose posterior expectation we seek to compute, and $c$ is a constant.
If the Stein equation has a solution, then $c$ \emph{must} be the value of the posterior expectation we seek.
This has motivated several efforts to numerically solve the Stein equation, as a more direct alternative to first approximating $p$ (for example using Markov chain Monte Carlo) and then using the approximation of $p$ to approximate the expectation of interest.
In this context kernel methods are typically used \citep{oates2017control,south2022semi} and in particular the kernel should have smoothness that is two orders higher than that of the function $f$ whose expectation is of interest, since the Stein equation is a second-order PDE.
The generalisation of the Stein equation to Riemannian manifolds was considered in \cite{barp2022riemann}, who advocated for the use of kernels with customisable smoothness that reproduce Sobolev spaces of functions on the manifold, such as the (manifold generalisation of the) \M model.
The connection between the \M model and Sobolev spaces is set out in Section \ref{SecNAaAT}.

\section{Flexible Modelling with Mat{\'e}rn} \label{secModeling}

One might object that the Mat{\'e}rn model is insufficiently flexible for many statistical applications, being limited to scalar-valued random fields that are stationary,  isotropic and Gaussian. 
However, the Mat{\'e}rn model is also an important building block for many more sophisticated models, some of which will now be described. 
This is a rich literature, and our discussion is necessarily succinct; an extended version of this section can be found in Appendix \ref{Appendix} of the Supplementary Material.

\subsection{Scalar Valued Random Fields} 

Let us start by discussing models for scalar-valued random fields that build on the \M model.
Note that one can trivially introduce non-zero mean functions into the \M model, or combine (additively or multiplicatively) kernels to obtain a potentially more expressive kernel; we will not dwell on either point.

To relax the isotropy assumption of the \M model, \cite{allard2016anisotropy} consider scale mixtures that take into account preferential directions in which spatial dependence develops.
On the other hand, the case of space-time models requires special treatment, and non-separable versions of the Mat{\'e}rn kernel are described in \citet{gneiting1,zast-porcu}.

The stationarity assumption was relaxed in a parametric manner in \citet{paciorek2006spatial}, and then in a nonparametric manner in \citet{roininen2019hyperpriors}.
An attempt to strike a balance between the computational tractability of parametric models and the flexibility of nonparametric models was reported in \citet{wilson2016deep}, who proposed \emph{input warping} to transform the inputs to the \M model using a neural network. 

The Gaussian assumption can be relaxed through \emph{output warping}, meaning transformation of the form  $\tilde{Z}(\bm{x}) = w(Z(\bm{x}))$ where $w(\cdot)$ is a nonlinear map from $\R^d$ to $\R^d$.
The covariance function of $\tilde{Z}$ will not be Mat\'{e}rn in general, when the covariance function of $Z$ is \M, but if $w$ is sufficiently regular then the smoothness properties of $Z$ transfer to $\tilde{Z}$.
The question of whether there exist non-Gaussian processes whose covariance function is nevertheless of Mat\'{e}rn class was answered positively in \citet{aaberg2011class}. 
 \citet{yan2018gaussian} have proposed \emph{trans-Gaussian} random fields with Mat{\'e}rn covariance function. \citet{bolin2014spatial} and subsequently \citet{wallin2015geostatistical} provided SPDE-based constructions for non-Gaussian Mat{\'e}rn fields.  
 General classes of non-Gaussian fields with covariance $g({\cal M}_{\nu,\alpha})$, for $g(\cdot)$ a suitable function that preserves  the positive definiteness and smoothness properties of the 
 Mat{\'e}rn model, have been provided for instance by \citet{Palacios:Steel:2006,Xua:Genton:2017,Bevilacqua_et_al:2021,mmaa}. 

An important extension of the \M model, which has received recent attention, is to random fields on spaces for which classical notions of smoothness are not well-defined.
For example, \citet{anderes2020isotropic} consider graphs with Euclidean edges, equipped with either the geodesic distance over the graph, or the resistance metric. 
\citet{menegatto2020gneiting} provide a generalisation of this setting by considering quasi-metric spaces.  \citet{https://doi.org/10.48550/arxiv.2205.06163} adopt a different approach to build random fields with their covariance structure on metric graphs. 
Space-time version of the Mat{\'e}rn model, for graphs with Euclidean edges, have been considered by \citet{tang2020space} and \citet{porcu2022nonseparable}. 
These efforts considerably extend the applicability of the \M model.


The Mat\'ern covariance function decays exponentially with distance, which can be inappropriate for modelling processes that involve long memory. 
Several approaches have been developed to modify the tails of the Mat{\'e}rn correlation function while preserving many of its desirable characteristics; we describe these in Section \ref{surr}.

\cite{guinness2022inverses} considers Gaussian random fields defined for lattices $\mathbb{Z}^d$ with a covariance function that is the restriction of the Mat{\'e}rn covariance to $\mathbb{Z}^d$. The resulting spectrum is smoothed version of the spectral density associated with the Mat{\'e}rn covariance. For this specific situation, the SPDE approximation can overestimate the scale, $\alpha$. Yet, it is not clear how this message extends to Gaussian fields that are continuously indexed in $\R^d$.


\subsection{Vector-Valued Random Fields}
 There has been a plethora of approaches related to multivariate spatial modeling, and the reader is referred to \citet{Genton:Kleiber:2014}. Here, the isotropic covariance function $\bm{K}: [0,\infty) \to \R^{p \times p}$ is matrix-valued. The elements on the diagonal, $K_{ii}$, are called \emph{auto-covariance} functions, and the elements $K_{ij}$, $i\ne j$, are called \emph{cross-covariance} functions. 
\citet{Gneiting:Kleibler:Schlather:2010} proposed a multivariate \M model 
\begin{equation}\label{stma}
K_{ij}(x) = \sigma_{ii} \sigma_{jj} \rho_{ij} {\cal M}_{\nu_{ij},\alpha_{ij}} (x), \qquad x \ge 0, \quad 
\end{equation}
where $\sigma_{ii}^2$ is the variance of $Z_i$, the $i$th component of a multivariate random field in $\R^p$, and $\rho_{ij}$ is the collocated correlation coefficient. There are restrictions on the parameters $\nu_{ij}, \alpha_{ij}$ and $\rho_{ij}$ required to ensure positive definiteness, and often the restrictions on the collocated correlations coefficients $\rho_{ij}$ are rather strict. 
This last remark has motivated alternative approaches, and the reader is referred to \citet{Apanasovich} and more recently to \citet{emery2022new}. Extensions to multivariate space-time Mat{\'e}rn structures have been provided by 
 \citet{allard2022fully} and through a technical approach by \citet{porcu2022criteria}. Multivariate nonstationary Mat{\'e}rn functions have been proposed by 
\citet{kleiber2012nonstationary}. 
Multivariate Mat{\'e}rn models with \emph{dimple} effect have been studied by \citet{alegria2021bivariate}; a `dimple' in a space-time covariance model refers to the case when $\text{Cov}( Z(\bm{x}, t), Z(\bm{x}', t') )$ is bigger than $\text{Cov}(( Z(\bm{x}, t), Z(\bm{x}', t))$, which requires special mathematical treatment. 

Multivariate Mat{\'e}rn modeling on graphs has been recently investigated in \citet{dey2022graphical}, who  propose a  class of multivariate graphical Gaussian processes through {\em stitching}, a construction that gets multivariate covariance functions from the graph, 
and ensures process-level conditional independence between variables. 
When coupled with the Matérn model, stitching yields a multivariate Gaussian process whose univariate components are Matérn Gaussian processes, and which agrees with process-level conditional independence as specified by the graphical model. Stitching can offer massive computational gains and parameter dimension reduction. An ingenious approach to Gaussian process construction involving the Mat{\'e}rn covariance function has been recently proposed  by \citet{li2021multi}, who considered a product space involving the $d$-dimensional Euclidean space cross an abstract set that allows to index group labels.

\subsection{Directions, Shapes and Curves}
The Mat{\'e}rn model has an important role in the study of directional processes, with \citet{banerjee_directional}  formalising the notions of directional \emph{finite difference processes} and \emph{directional derivative processes} with special emphasis on the Mat{\'e}rn model. 
The \M model also has a role in shape analysis, where \citet{banerjee2006bayesian} introduced \emph{Bayesian wombling} to measure \emph{spatial} gradients related to curves through `wombling' boundaries, and approach taken further in \citet{halder2023bayesian}. 
The smoothness properties of the Mat{\'e}rn model are ideally suited to such a framework.
Modeling approaches to \emph{temporal} gradients using the Mat{\'e}rn model have been proposed by \citet{quick2013modeling}.  
Related to these approaches, the smoothness parameter $\nu$ of the \M model plays a central role in the recent paper by \citet{halder2023bayesian}, who analyse random surfaces in order to explain latent dependence within a response variable of interest. 


\smallskip

This represents a short tour of \emph{statistical} applications of the \M model, but its reach goes well beyond statistics, and we explore the importance of the \M model to related fields next.

\section{The Mat{\'e}rn Model Outside Statistics} \label{sec5}

This section explores the impact of the \M model on numerical analysis and approximation theory (Section \ref{SecNAaAT}), machine learning (Section \ref{SecML}), and probability theory (Section \ref{PTaSP}).

\subsection{Numerical Analysis and Approximation Theory}\label{SecNAaAT}

The problem considered here is to {\em reconstruct} a
  real-valued function $f$ defined on a domain 
  $D \subset \R^d$ from given {\em data values}
  $y_i = f(\bm{x}_i)$  available at a set
  $X_n =\{\bm{x}_1,\dots,\bm{x}_n\}$
  of distinct {\em data locations}. 
  In contrast to the statistical exposition in Section \ref{seclik}, from a numerical analysis standpoint these data are not assumed to be random in any way. 
Nevertheless, many of the mathematical expressions that we previously motivated from a statistical perspective appear also in the solution of this numerical task.
The data vector $\boldsymbol{Z}_n$ is reinterpreted as 
$\boldsymbol{Z}_n=(f({\bm x}_1),\ldots,f({\bm x}_n))^{\top}$
and the task is to approximate the value $f(\bm{x})$ of 
the unknown function $f$ at an unsampled location $\boldsymbol{x}\in D \setminus X_n$. 
A natural solution is a minimal-norm interpolant
\begin{align*}
    s_{f,X_n,K} = \argmin_{s \in \mathcal{H}(K)} \|s\|_{\mathcal{H}(K)} \quad \text{s.t.} \quad \begin{array}{l} s(\bm{x}_i) = f(\bm{x}_i), \\ \hspace{30pt} i = 1,\dots,n , \end{array} 
\end{align*}
which we recall was the third optimality property referred in Section \ref{Sec3}.
Thus, using again the {\em kernel matrix} 
$\bR_{n}=[K(\bm{x}_i,\bm{x}_j)]_{i,j=1}^n$,
the system $\bR_{n}\boldsymbol{c}_n=\bZ_n$ is solved for a 
fixed coefficient vector $\boldsymbol{c}_n$ that determines a
linear combination
$$
s_{f,{X_n},K}(\bx)=\sum_{i=1}^n c_iK(\bm{x}_i,\bx), \qquad \bx \in D,
$$
in the span of the {\em translates}
$K(\bm{x}_i,\cdot)$. 
This follows
easily from the reproduction formula (\ref{eqrepro}) 
and (\ref{eqHKKK}).
The above formula is identical to (\ref{blup}) when setting $\bx=\bx_0$, and
the resulting value 
$s_{f,{X_n},K}(\bx)$ is interpreted as a numerical approximation to $f(\bx)$.
The log-likelihood function (\ref{eq:17}) can equivalently be viewed as penalising the norm of the interpolant, since $\|s_{f,{X_n},K}\|^2_{{\cal H}(K)}=\boldsymbol{Z}_n^{\top} \bR_{n}^{-1}\boldsymbol{Z}_n.$

The fourth optimality principle in
Section \ref{Sec3} corresponds here to the fact that the
norm of the error functional 
$\epsilon_{\bx}\;:\;f\mapsto f(\bx)-s_{f,{X_n},K}(\bx)$
in the dual space $\mathcal{H}(K)^*$ of ${\cal H}(K)$ is minimal under all
linear reconstruction algorithms in ${\cal H}(K)$
that use the same data $\bm{Z}_n$.  
The key tool is 
the {\em power function} $P_{K,{X_n}}$, 
defined for all $\bx \in D$ by 
\begin{align*}
& \hspace{-5pt} P_{K,{X_n}}(\bm{x}) \\
& = \sup\left\{f(\bm{x}) \;:\; f\in {\cal H}(K),\;f({X_n})=0,\;\|f\|_{{\cal H}(K)}\le 1
\right\}
\end{align*}
It has the property $P_{K,{X_n}}(\bm{x})
=\|\epsilon_{\bx}\|_{{\cal H}^*(K)}$ 
and leads to optimal
error bounds of the form
$$ \left|f(\bm{x}) - s_{f,{X_n},K}(\bm{x}) \right| \le P_{K,{X_n}}(\bm{x}) \|f\|_{{\cal H}_K} .
$$
for all $\bx\in D$ and $f\in {\cal H}(K)$. 
It can be numerically calculated
using the kernel matrix based on ${X_n}\cup \{\bm{x}\}$, but we omit the detail.
Strikingly, the power function
  coincides with the square root of the {\em kriging variance}
  \citep{scheuerer2013interpolation},
  giving the variance of the kriging error at $\bm{x}$ for given data locations
  ${X_n}$ and kernel $K$.


Analysis of the approximation error in this context thus reduces to analysis of the power function, and in turn analysis of the space $\mathcal{H}(K)$.
From (\ref{eqHK}) and (\ref{stein1}), the RKHS generated by the \M 
kernel ${\cal M}_{\nu,1}$ has the inner product 
\begin{equation}\label{eqHKM}
    \langle f,g\rangle_{{\cal H}({\cal M}_{\nu,1})}=\int_{\R^d}\frac{\hat f(\bm{\omega})
\overline{\hat g(\bm{\omega})}}{(1+\|\bm\omega\|^2)^{\nu+d/2}}{\rm d}\bm{\omega}
\end{equation}
up to constants, 
which we recognise as the inner product of the classical \emph{Sobolev space} $W_2^{\nu+d/2}(\R^d)$.
By the Sobolev embedding theorem, the elements of this space are well-defined continuous functions whenever $\nu > 0$.
This space is a canonical setting for mathematical analysis of PDEs, a connection that we trailed in Section \ref{SecQua}.
Summarising, the use of \M kernels
  yields optimal recovery techniques for functions in Sobolev spaces
  from given sampled data. 
  Generalised recoveries using derivative data
produce \emph{meshless} numerical methods for solving PDEs
in Sobolev spaces, including the {\em symmetric collocation} method 
which uses derivative data for the PDE based on \citet{wu:1992-1}, and shares similar Hilbert space optimality properties \citet{schaback:2015-3}.
The use of the \M kernel is strongly motivated by the fact that PDE theory
often implies that solutions lie in Sobolev spaces. 
On the other hand, there are also good arguments to replace \M kernels by polyharmonics \cite{schaback:2016-1,davydov-schaback:2019-1}.
 
Plenty of other results on deterministic
recovery problems using kernels can be found in \citet{wendland2004scattered},
while applications are in \citet{schaback2006kernel} and 
MATLAB programs combined with the essential theory are in \citet{fasshauer-mccourt:2015-1}.

In numerical analysis and approximation theory, Mat\'ern and other kernels are normally used for rather large values of their smoothness parameter, because they seek to solve an interpolation rather than a regression task. 
\citet{narcowich2006sobolev} proved that convergence rates then depend on the minimum of the smoothness of the function $f$ providing the data and the kernel; 
a \emph{misspecified} \M kernel, for which the smoothness parameter $\nu$ is taken to be too large relative to the smoothness of $f$, produces an error that converges at the same rate as we would have achieved had $\nu$ been correctly specified. 
On the other hand, \citet{tuo2020kriging} prove in the same setting that the prediction error becomes more sensitive to the space-filling property of the design points. In particular, optimal convergence rates require also that the \emph{quasi-uniformity} of the experimental design is controlled.

Of course, the use of kernels
in numerical analysis and approximation theory requires estimation
of kernel parameters. The quantity $\sigma$ does not arise
in the correlation matrix $\bR_{n}$, but the scale parameter $\alpha$ has
a strong influence on the error of the interpolant. 
There is a vast literature on \emph{scale estimation} that partially builds on
statistical notions like ML (see references in Section 3).
On the other hand, specific alternatives to the \M model, such as the polyharmonic kernels of Section \ref{SecPoly}, are able to bypass scale estimation due to the remarkable property that the interpolant is independent of the value of the scale parameter used.
See \citet{wendland2004scattered} and Section \ref{SecPoly}.

\subsection{Machine Learning}\label{SecML}

Kernel methods are a major strand of machine learning research, where kernels are routinely used to solve a variety of supervised and unsupervised learning tasks.
Compared to the interpolatory setting of Section \ref{SecNAaAT}, data in machine learning are usually observed with noise, necessitating either a likelihood or a loss function to be specified.

The \M model is often convenient for the analysis of kernel methods; for example, \citet{tuo2020improved} provide sufficient conditions for the rates of convergence of the Mat{\'e}rn kernel ridge regression to  exceed the standard minimax rates under both the $L_2$ norm and the norm of the RKHS.
However, the presence of noise in the data can pose a substantial challenge to selection of smoothness parameters such as $\nu$ in the \M model.
\citet{karvonen2022asymptotic} proves that the ML estimate of $\nu$ cannot asymptotically {\em undersmooth} the truth under fixed domain asymptotics; that is, if the true regression function has a Sobolev smoothness $\nu_0 + d/2$, then the smoothness parameter estimate cannot be asymptotically less than $\nu_0 + d/2$, but this in itself it not compelling motivation to use ML \citep{karvonen2022maximum}. 
As a result of these additional challenges, standard practice is to keep the kernel general as far as possible when developing methodology, and as far as possible to learn a suitable form for the kernel using the data and model selection criteria.  
However, recent machine learning methodology for non-Euclidean data hinges on the SPDE approach, and as a consequence the Mat\'{e}rn and related models are explicitly being used.
As the types of data that researchers seek to analyse become more heterogeneous and structured, there has been a demand for flexible Gaussian process models defined on such non-Euclidean domains as manifolds and discrete, graph-based domains.
Under the framework of Gaussian processes, \citet{borovitskiy2020matern} 
proposed to avoid numerical solution of the SPDE \eqref{SPDE} and instead to work with a finite-rank approximation to the Gaussian process model.
Specifically, they consider the SPDE in (\ref{SPDE}) appropriately adapted to a Riemannian manifold $M$, for which the corresponding Mat{\'e}rn model admits a series expansion of the type
$$ \sum_{n=0}^{\infty} \left ( \frac{2 \nu}{\alpha^2} + \lambda_n \right )^{-\nu-d/2} f_{n}(\bm{x})f_n(\bm{x}'), \qquad \bm{x},\bm{x}' \in M $$
where $\{\lambda_n\}_{n=0}^{\infty}$ and $\{f_n\}_{n=0}^{\infty}$ are, respectively, the sequences of eigenvalues and eigenfunctions from the Laplace--Beltrami operator $- \Delta_M$.
The authors propose to first solve numerically for the leading eigenfunctions $\{f_i\}_{i=0}^n$ of the Laplace--Beltrami operator, and then working with a finite-rank Gaussian process whose realisations are linear combinations of the $\{f_i\}_{i=0}^n$.
Though solving the eigenproblem may be harder than numerically solving the SPDE, the authors argue that caching of the eigenfunctions can lead to a cost saving in settings where multiple tasks are to be solved on the same manifold.
Such an approach is ingeniously extended to undirected graphs by \citet{pmlr-v130-borovitskiy21a}, and has had  a direct impact on Gaussian processes defined on neural networks \citep{jensen2020manifold}, pathwise conditioning of Gaussian processes \citep{wilson2021pathwise}, simulation intelligence in AI \citep{lavin2021simulation} and extension to kernel methods withing graphs cross time \citep{nikitin2022non}. Other applications include Thomson sampling in neural information systems \citep{vakili2021scalable}, Bayesian optimisation in robotics \citep{jaquier2022geometry}, and Gaussian processes regression on metric spaces \citep{koepernik2021consistency}.

\subsection{Probability Theory and Stochastic Processes}\label{PTaSP}
The Mat{\'e}rn model is well-studied from a probability theory and stochastic process viewpoint. 
From the perspective of regularity, \citet{scheuerer2010regularity} summarises the properties of Gaussian random fields with Mat{\'e}rn covariance functions; sample paths are $k$-times differentiable in the mean-square sense if and only if $\nu >k$. Under the same condition, the sample paths have (local) Sobolev space exponent being identically equal to $k$. 
Further, a Gaussian random field with Mat{\'e}rn covariance has fractal dimension that is identically equal to $\min (\nu, d)$, for $d$ being the dimension of the Euclidean space on which the random field is defined. 
For non-Gaussian random fields with Mat{\'e}rn covariance, continuity properties are studied by \citet{kent1989continuity}.

Several other properties of the \M model have been investigated.
\citet{kelbert2005fractional} study fractional random fields under the scenario of stochastic fractional heat equations under a Mat{\'e}rn model; see also \citet{leonenko2011fractional}.
Random fields defined on the unit ball embedded in $\R^d$, with a covariance function that is the restriction of the Mat{\'e}rn model to a finite range, were studied in \citet{leonenko2022spectral}. 
Tensor-valued random fields with an equivalent class of Mat{\'e}rn covariance functions were studied in \citet{leonenko2017matern}.
\citet{terdik2015angular} considers angular spectra for non-Gaussian random fields with \M covariance function. 
A recent contribution \citep{terdik2022spatiotemporal} provides interesting connection between the Mat{\'e}rn model and certain Laplacian ARMA representations of a class of stochastic processes. \citet{lilly2017fractional} show that the Matérn process is a damped version of fractional Brownian motion. 
\citet{lim2009generalized} study random fields with a generalised Mat{\'e}rn covariance obtained as the solution to the fractional stochastic differential equation with two fractional orders, enabling the authors to 
deduce the sample path properties of the associated random field. 
Space-time extensions of Mat{\'e}rn random fields through stochastic Helmholtz 
equations are provided by \citet{angulo2008spatiotemporal}. 

According to N. Leonenko\footnote{Personal Communication, January 2023.}, a major contributor to this literature, the importance of Mat{\'e}rn model is based on the Duality theorem \citep[][Theorem 1]{harrar2006duality} which provides an explicit relation between certain classes of characteristic functions of symmetric random vectors and their density. Specifically, the spectral density associated with the Mat{\'e}rn model is by itself a covariance function, called the Cauchy 
or inverse multiquadric covariance function, that allows to parameterise the Hurst effect of the associated Gaussian random field. 


\smallskip

This completes our tour across the scientific landscape through the lens of the \M model.
Our attention turns now to the future, and promising enhancements that can be made to the \M model.

\section{Enhancements of the \M Model}\label{surr}

 \label{sec6}



This section described \emph{enhancements} of the \M model; covariance functions that share (at least partially) the local properties of the Mat{\'e}rn model while providing additional features and functionality. 
Here we first introduces the models one at a time, with critical commentary on their features deferred to Section \ref{sec7}.

\subsection{Models with Compact Support} \label{sec-comp-sup}


Compactly supported covariance models have a long history that can be traced back to \citet{Askey:1973}, who proposed  the kernel 
\begin{equation}
    \label{askey}
    {\cal A}_{\mu,\beta}(x) = \left ( 1- \frac{x}{\beta}\right )_+^{\mu}, \qquad x \ge 0,
\end{equation}
with $\beta$ and $\mu$ being strictly positive, and where $(x)_+=\max(0,x)$ is the {\em truncated power}. 
It was shown in that work that ${\cal A}_{\mu,\beta}$ belongs to $\Phi_d$ for all $\beta>0$ if and only if $\mu \ge (d+1)/2$. 
Clearly, the mapping $\bm{x} \mapsto {\cal A}_{\mu,\beta}(\|\bm{x}\|)$ is compactly supported over a ball with radius $\beta$ embedded in $\R^d$. 
As a result, covariance matrices contain exact zero entries whenever the associated states $\bm{x}_i$ and $\bm{x}_j$ satisfy $\|\bm{x}_i - \bm{x}_j\| \geq \beta$; the computational advantages of this sparsity are discussed further in Section \ref{subsec: big data and comp}.

Matheron's \emph{mont{\'e}e} and \emph{descente} \citep{matheron65} 
approach was applied by \citet{Wendland:1995} to the Askey functions, obtaining compactly supported covariance functions with higher-order smoothness that are truncated polynomials as functions of $\|\bm{x}\|$.
This strategy was unable to generate integer-order Sobolev spaces in even space dimensions, a problem that was resolved in \citet{schaback2011missing} who identified the `missing' Wendland functions. 
A unified view of Wendland functions was provided by \citet{gneiting1}. \citet{zast2002} provided necessary and sufficient conditions for a general class encompassing both ordinary and missing Wendland functions. \citet{buhmann2001new} provided a generalisation of Wendland functions, with sufficient parametric conditions that allow the new class to belong to $\Phi_d$ for a given $d$. Those functions, termed \emph{Buhmann} functions, were then studied by \citet{zastavnyi2006some} and subsequently by \citet{zastavnyi2017positive,porcu2016buhmann} and \citet{faouzi2020zastavnyi}. 
Alternative representations and properties of the Wendland functions have been studied by \citet{hubbert2012closed} and \citet{chernih2014closed}.
Extensions of the Wendland functions to multivariate \cite{porcu2013radial,daley2015classes}, spatio-temporal \cite{porcu2020nonseparable} and non-stationary processes \cite{kleiber2015nonstationary} have also been developed.


A more technical discussion follows, in which we introduce two further classes of correlation functions with compact support, each of which will be the subject of discussion in Section \ref{sec7}.
\begin{enumerate}
\item The generalized Wendland (${\cal GW}$) family \citep{Gneiting:2002b,zastavnyi2006some} contains correlation functions with compact support that, as in the Mat\'ern model, admit a continuous parameterisation of smoothness of the underlying Gaussian random field. 
The ${\cal GW}_{\kappa,\mu,\beta}$ model depends on parameters $\kappa \ge 0$ and $\mu,\beta>0$ through the identity
\begin{multline} \label{WG4*}
	{\cal GW}_{\kappa,\mu, \beta}(x)\\= 
	  \frac{\Gamma(\kappa)\Gamma(2\kappa+\mu+1)}{\Gamma(2\kappa)\Gamma(\kappa+\mu+1)2^{\mu+1}}  {\cal A}_{\kappa+\mu,\beta^2} \left( x^2 \right) \\
	   \times  {}_2F_1\left(\frac{\mu}{2},\frac{\mu+1}{2};\kappa+\mu+1;{\cal A}_{1,\beta^2} \left( x^2 \right) \right),
	\end{multline}
 where 
$\mu\geq  (d+1)/2+\kappa$ is needed for ${\cal GW}_{\kappa,\mu, \beta}$ to belong to the class $\Phi_d$
and  ${}_2F_1(a,b,c,\cdot)$ is the Gaussian hypergeometric function \citep{Abra:Steg:70}.
Sample paths of the ${\cal GW}_{\kappa,\mu, \beta}$ model are $k$ times mean-square differentiable, in any direction, if and only if $\kappa>k-1/2$ \citep{Gneiting:2002b}, so that $\kappa$ plays the role of the smoothness parameter in this model.
When $\kappa=k \in \mathbb{N}$,
${\cal GW}_{k,\mu, \beta}$
 factors into the product of the Askey function ${\cal A}_{\mu+k,\beta}$ with a polynomial of degree $k$.  
This model includes the Wendland functions ($\kappa=k$, a positive integer), as well as the  missing Wendland functions ($\kappa=k+1/2$). 
Theorem 1(3) in \citet{BFFP} implies that the RKHS induced by ${\cal GW}_{\kappa/2-(d+1)/4,\mu,\beta}$, with $\kappa \ge (d+1/2)$, is norm-equivalent to the Sobolev space $W_{2}^{\kappa}(\mathbb{R}^d)$.


\item The Gauss hypergeometric ($\mathcal{GH}$) family  \citep{emery2022gauss} is  defined as
\begin{multline}
    \label{emery}
    \mathcal{GH}_{\kappa,\delta,\gamma,\beta}(x) \\ = \frac{\Gamma(\delta-d/2)\Gamma(\gamma-d/2)}{\Gamma(\delta-\kappa+\gamma - d/2) \Gamma(\kappa-d/2)} \\ \times  {\cal A}_{\delta-\kappa + \gamma-d/2+1,\beta^2} \left( x^2 \right)\times\\
    {}_2F_1\left( \delta-\kappa; \gamma-\kappa; \delta-\kappa + \gamma -d/2; 
    {\cal A}_{1,\beta^2} (x^2) 
    \right).
\end{multline}
 This model has four parameters and it belongs to the class $\Phi_d$ for every positive $\beta$ provided  $\kappa>d/2$ with 
  $$ 2(\delta-\kappa)(\gamma-\kappa) \geq \kappa, \quad \text{and} \quad
2(\delta+\gamma)\geq 6\kappa+1. $$
 Sample paths of the $\mathcal{GH}_{\kappa,\delta,\gamma,\beta}$ model are  $\lceil k/2\rceil$ times mean-square differentiable, in any direction, if and only if $\kappa>(k+d)/2$.
 The parameter $\kappa$ thus also controls the smoothness of samples from this model. 


\end{enumerate} 

\noindent The importance of the $\mathcal{GW}$ and $\mathcal{GH}$ models is discussed in Section \ref{sec7}.

\subsection{Models with Polynomial Decay}

Correlation models with polynomial decay such as the 
generalized
Cauchy \citep{gneiting2004stochastic} or the Dagum models \citep{berg2008dagum}
can be useful when modelling data with long-range dependence.
However, in using these correlation models
one loses control over the   differentiability of the  the sample paths,  a key property of the Mat\'ern model.
\citet{ma2022beyond} recently proposed a modification of the Mat\'ern class that allows for polynomial decay, while maintaining the local properties of the conventional Mat\'ern model.  The correlation function associated to this model is given by 
\begin{equation}
\label{anciano}
    \mathcal{CH}_{\nu,\eta,\beta}(x) = \frac{\Gamma(\nu+\eta)}{\Gamma(\nu)} \, \mathcal{U}\left( \eta,1-\nu, \nu \left(\frac{x}{\beta} \right)^2 \right), \; x\geq 0, 
\end{equation}
where $\mathcal{U}$ is the confluent hypergeometric function of the second kind \citep{Abra:Steg:70}.
Here $\nu>0$ controls mean-square differentiability near the origin, as in the Mat\'ern case, while $\eta>0$ controls the heaviness of the tail.
The construction (\ref{anciano}) is based on a scale mixture of (a reparameterised version of) the Mat\'ern model 
involving the inverse-gamma distribution. 
 \citet{ma2022beyond} have shown that this class is particularly useful for extrapolation problems where large distances are predominant.

\subsection{Polyharmonic Kernels}\label{SecPoly}

Our catalogue of enhancements of the Mat{\'e}rn model finishes with
{\em polyharmonic kernels}, defined as 
\begin{equation}\label{eqKmd}
H_{\nu,d}(x):=
\left\{
\begin{array}{ll}
x^{2\nu-d}\log x &\hbox{ for } 2\nu-d \in 2\Z\\ 
x^{2\nu-d}& \hbox{ else }
\end{array} 
\right\}
\end{equation}
up to the sign $(-1)^{\lfloor \nu-d/2 \rfloor +1}$. 
As a function of $x=\|\bm{x}\|$, $\bm{x} \in \R^d$,
the Mat\'ern kernel ${\cal M}_{\nu-d/2,1}$ starts with even powers of $x$ followed by $H_{\nu,d}$, and in this sense the two models are related.
Up to a constant factor, the generalised Fourier transform of $H_{\nu,d}(\|x\|)$ on $\R^d$ 
is $\|\omega\|^{-2\nu}$, and then a scale parameter is just another constant factor. This makes kernel-based interpolation by
polyharmonics scale-independent. Compare with (\ref{stein1}) to see the connection to ${\cal M}_{\nu-d/2,\alpha}$ in Fourier space.
\citet{stein2004equivalence} provides a formal connection between polyharmonic kernels, for which the name {\em power law} covariance functions is also used, and the Mat{\'e}rn model.
%
Polyharmonic kernels
are {\em conditionally positive definite} of order 
$\lfloor \nu -d/2 \rfloor+1$; for a technical definition see \citet{wendland2004scattered}.
Instead of Hilbert Spaces, polyharmonic kernels generate
{\em Beppo--Levi spaces}, which share similarities to Sobolev spaces modulo that an additional polynomial space has to be added to enable prediction (Section \ref{Sec3}) and interpolation (Section \ref{SecNAaAT});
see \citet{wendland2004scattered}.
In general, polyharmonic kernels arise as covariances in {\em 
fractional} Gaussian fields, including forms of Brownian motion \citep[][Theorem 3.3]{lodhia-et-al:2016}. 

\smallskip

Next our attention turns to a critical discussion of whether such enhancements to the \M model are needed.





\section{Are Enhancements of the Mat{\'e}rn Model Useful?}
\label{sec7}

This final section provides critical commentary on the \M model and the enhanced versions of the model introduced in Section \ref{surr}.

\subsection{Rigorous Generalisation of the Mat{\'e}rn Model }

The Mat{\'e}rn model does not allow for compact support, hole effects (oscillations between positive and negative values) at large distances, or slowly decaying tails suitable for modeling long-range dependence. 
 Most of the enhancements in Section \ref{surr} aim to resolve these kind of issues; here we describe how the $\mathcal{GW}$, $\mathcal{GH}$ and $\mathcal{CH}$ models can be viewed as rigorous generalisations of the Mat{\'e}rn model.

 \citet{bevilacqua2022unifying} have shown that the
Mat{\'e}rn model is a limit case of a rescaled version of  the ${\cal GW}$ model.
In particular they have considered the model  $\widetilde{{\cal GW}}$ defined as
$$\widetilde{{\cal GW}}_{\kappa,\mu,\beta}(x)  = {\cal GW}_{\kappa,\mu,\beta \left(\frac{\Gamma(\mu+2\kappa +1)} {\Gamma(\mu)}\right)^{\frac{1}{1+2\kappa}}}(x), \qquad x \ge 0, $$
and proved that 
 $$\lim_{\mu\to\infty} \widetilde{{\cal GW}}_{\kappa,\mu,\beta}(x)={\cal M}_{\kappa+1/2,\beta}(x),\quad \kappa\geq 0,$$
with uniform convergence over the set $x\in (0,\infty)$.
The parameter $\mu$  thus allows for switching from compactly supported to globally supported models, and can either be fixed to ensure
sparse correlation matrices, or can be estimated based on the dataset.
However, this equivalence applies only to smoothness parameters greater than or equal to $1/2$ in the \M model, so the full range of the smoothness parameter is not covered.
This is unfortunate, since the fractal dimension
 \citep[a widely used measure of roughness  of the sample paths for time series and spatial data;][]{gneper2012} is fully parameterised using the Mat{\'e}rn model when the smoothness parameter lies between $0$ and $1$.
As a consequence,  the ${\cal GW}$ (or $\widetilde{{\cal GW}}$) model   cannot fully parameterise the fractal dimension of the random field. This kind of issue can be solved with the 
$\mathcal{GH}$ model, which includes the ${\cal GW}$ model as a special case \cite{emery2022gauss}:
$$ \mathcal{GH}_{\frac{d+1}{2}+\nu,\frac{d+\mu+1}{2}+\nu,\frac{d+\mu}{2}+1+\nu,\beta}(x)= {\cal GW}_{\nu,\mu,\beta}(x)$$
Letting $\beta$, $\delta$ and $\gamma$ tend to infinity in such a way that $\beta/\sqrt{4\delta\gamma}$ tends to $\alpha>0$, the $\mathcal{GH}$ model (\ref{emery}) converges uniformly to the Mat\'ern model $\mathcal{M}_{\kappa-d/2, \alpha}(x)$, and in this case the \emph{full range} of the smoothness parameter of the Mat\'ern model is covered. 

The  Mat\'ern model also arises as a special limit case of the $\mathcal{CH}$ model.
Specifically, \citet{ma2022beyond} show that 
 $$\lim_{\eta\to\infty} \mathcal{CH}_{\nu,\eta,2\sqrt{\nu(\eta+1)}\beta}(x)={\cal M}_{\nu, \beta}(x),$$
 with convergence being uniform on any compact set.

 The \emph{turning band} operator of \citet{matheron63} can be applied to a correlation function to create hole effects
 while retaining positive definiteness of the kernel. 
 An argument in Schoenberg proves that, for an isotropic correlation in $\R^d$, the correlation values cannot be smaller than $-1/d$ \citep{schoenberg}.
 Since the Mat{\'e}rn model is a valid model for all $d$, this implies that the application of turning bands to the Mat{\'e}rn model will not provide any hole effect. On the other hand, the ${\cal GW}$ and ${\cal GH}$ models allow for such an effect.

\subsection{Estimation of Enhanced Models} \label{estim}

ML estimation for the Mat{\'e}rn model are well-understood; here we discuss the extent to which similar results can be obtained for enhancements of the \M model.

In the context of increasing domain asymptotics, parameters of the ${\cal GW}$ and ${\cal CH}$ models can be estimated consistently using ML and  the associated asymptotic distribution is known; see Section \ref{ida}.

In the context of fixed domain asymptotics, similar to the classical Mat{\'e}rn model, the parameters of the these enhanced models cannot be consistently estimated.
For instance, 
\citet{BFFP} show that the microergodic parameter of the   covariance model
$\sigma^2{\cal GW}_{\kappa,\mu,\beta}$, assuming $\kappa$ and $\mu$ known,
is given by $ {\rm micro}_{{\cal GW}}= {\sigma^2}/{ \beta^{2 \kappa+1}}.$
In addition they prove that 
for a zero mean Gaussian field  defined on  a  bounded infinite set $D\subset \R^d$ ($d=1, 2, 3$), 
with covariance model $\sigma^2_0{\cal GW}_{\kappa,\mu,\beta_0}$,
the ML estimator $\hat{\sigma}_{n}^{2}/\hat{\beta}_{n}^{2\kappa+1}$ of the microergodic parameter 
is strongly  consistent,  $i.e.$,  $$\hat{\sigma}_{n}^{2}/\hat{\beta}_{n}^{2\kappa+1}\stackrel{a.s.}{\longrightarrow} \sigma_{0}^{2}/\beta_{0}^{2\kappa+1}. $$ Additionally, 
for $\mu> (d+1)/2 +\kappa+ 3$,
its  asymptotic distribution is given by 
$$\sqrt{n}(\hat{\sigma}_{n}^{2}/\hat{\beta}_{n}^{2\kappa+1}-\sigma_{0}^{2}/\beta_{0}^{2\kappa+1})\stackrel{d}{\longrightarrow} \mathcal{N}(0,2(\sigma_{0}^{2}/\beta_{0}^{2\kappa+1})^2).$$
Analogous for the $\mathcal{GH}$ model proposed are not available at present.

Similarly, \citet{ma2022beyond} show that the microergodic parameter of the   covariance model
$\sigma^2{\cal CH}_{\nu,\eta,\beta}$, assuming $\nu$  known,
is given by $$ {\rm micro}_{{\cal CH}}= (\sigma^2 \Gamma(\nu+\eta))/ (\beta^{2 \nu}\Gamma(\eta)).$$
In addition they prove that 
for a zero mean Gaussian field  defined on  a  bounded infinite set $D\subset \R^d$ ($d=1, 2, 3$), 
with covariance model $\sigma^2_0{\cal CH}_{\nu,\eta_0,\beta_0}$,
the ML estimator $(\hat{\sigma}_{n}^{2}/\hat{\beta}_{n}^{2\nu})(\Gamma(\nu+\hat{\eta}_n)/\Gamma(\hat{\eta}_n))$ of the microergodic parameter 
is strongly  consistent,  $i.e.$,  $$\frac{\hat{\sigma}_{n}^{2} (\Gamma(\nu+\hat{\eta}_n) }{ \hat{\beta}_{n}^{2\nu} \Gamma(\hat{\eta}_n)}\stackrel{a.s.}{\longrightarrow} \frac{\sigma_{0}^{2}\Gamma(\nu+\eta_0)}{\beta_{0}^{2\nu} \Gamma(\eta_0)}$$ and,
if $\eta_0> d/2$,
its  asymptotic distribution is given by 
\begin{align*}
& \hspace{-40pt} \frac{\hat{\sigma}_{n}^{2} (\Gamma(\nu+\hat{\eta}_n) }{ \hat{\beta}_{n}^{2\nu} \Gamma(\hat{\eta}_n)}-\ \frac{\sigma_{0}^{2}\Gamma(\nu+\eta_0)}{\beta_{0}^{2\nu} \Gamma(\eta_0)} \\
& \stackrel{d}{\longrightarrow} \mathcal{N}\left(0,2\left( \frac{\sigma_{0}^{2}\Gamma(\nu+\eta_0)}{\beta_{0}^{2\nu} \Gamma(\eta_0)}\right)^2 \right).
\end{align*}
These results broadly support the use of ML plug-in estimates for these enhanced versions of the \M model; the issue of predictive performance is discussed next.




\subsection{Prediction with Enhanced Models}

If two Gaussian measures are equivalent then the associated predictions and mean squared errors are asymptotically identical (c.f. Section \ref{pprr}). 
To this end, recent results have sought to establish equivalence between Gaussian measures for the Mat{\'e}rn model and enhancements of the \M model.
\citet{BFFP} consider the  $\sigma_1^2{\cal GW}_{\kappa,\mu,\beta}$ model and 
show that 
for given $\sigma_1 \geq 0$, $\nu \geq 1/2$, and $\kappa \geq 0$,  if $\nu=\kappa+1/2$, $\mu > d+ \kappa+1/2$ and  
\begin{equation}\label{true}
\sigma_{0}^{2}\alpha^{-2\nu}=\left( \frac{\Gamma(2\kappa+\mu+1)}{\Gamma(\mu)} \right) \sigma_{1}^{2}\beta^{-(1+2\kappa)},
\end{equation}
then 
$P(\sigma^2_0{\cal M}_{\nu,\alpha})$ is equivalent to $P(\sigma^2_1{\cal GW}_{\kappa,\mu,\beta})$, for $d=1,2,3$, on the paths of $Z(\bm{x})$ for $\bm{x} \in D \subset \R^d$. 
Thus predictions made using the $\mathcal{GW}$ model with compact support are asymptotically identical to those made using the \M model.
Likewise, \citet{ma2022beyond} show that 
for a  given $\eta\geq d/2$ and $\nu \geq 0$,  if
\begin{equation}\label{true2} 
\sigma_{0}^{2}\alpha^{-2\nu}=\left( \frac{\Gamma(\nu+\eta)}{\Gamma(\eta)} \right) \sigma_{1}^{2}\left(\frac{\beta^2}{2}\right)^{-\nu},
\end{equation}
then 
$P(\sigma^2_0{\cal M}_{\nu,\alpha})$ is equivalent to $P(\sigma^2_1{\cal CH}_{\nu,\eta,\beta})$, for $d=1,2,3$, on the paths of $Z(\bm{x})$ for $\bm{x} \in D \subset \R^d$. 
Thus predictions made using the $\mathcal{GW}$ model with polynomial tail decay are asymptotically identical to those made using the \M model.




If interest is in the predictor \eqref{blup}, but not the predictive uncertainty resulting from the associated Gaussian random field, then it is interesting to note that the stationarity assumption of the \M model may not be needed.
\citet{stein2004} showed that, under suitable parametric conditions, one can 
consider $\alpha=0$ in the Mat{\'e}rn model, and this is equivalent to prediction using  the polyharmonic kernels $H_{\nu,d}$ in (\ref{eqKmd}). 
Theorem 1 in that work shows that if $d\le 3 $ and the parameter $\nu$ satisfies condition (2) therein (or $d = 1$), then {\em it is impossible to distinguish $\alpha>0$ from $\alpha=0$ on a bounded domain.} 
The above observation reflects the fact that prediction using polyharmonic kernels, like in Section \ref{SecNAaAT}, is scale-independent. 
This follows from homogeneity of the Fourier transform and eliminates the need for scale estimation in this context.

\subsection{Screening with Enhanced Models}

The screening effect extends also to enhanced versions of the \M model.
For regular schemes, Theorem 1 in \citet{porcu2020stein} shows that the ${\cal GW}$ model allows for an asymptotic screening effect when $\mu > (d+1)/2 +\kappa $.  
This condition is not restrictive, since  $\mu \ge (d+1)/2 +\nu $ is already required for ${\cal GW}_{\kappa,\mu,\beta}$ to belong to the class $\Phi_d$.
For irregular schemes the situations is more complicated. 
For example, for non-differentiable fields in $d=1$, Theorem $1$ in \citet{stein2} in concert with Theorem 1 in \citet{porcu2020stein} explains that the Askey model ${\cal GW}_{0,\mu,\beta}$ allows for a screening effect provided $\mu >1$. 
For $d=2$, Theorem $2$ in \citet{stein2} implies that the Askey model allows for screening provided that $\mu>3/2$.
The $\mathcal{GW}$ model satisfies Stein's condition in (1.3) of \citet{porcu2020stein}, which in turn allows the Stein hypothesis \eqref{condition1} to be verified.  

The numerical experiments in \citet{porcu2020stein} suggest that the screening effect is even stronger under enhanced models with compact support, compared to the standard Mat{\'e}rn model. 
This can deliver computational advantages, which we discuss next.


\subsection{Scalable Computation}
\label{subsec: big data and comp}






The eternal fight between statistical accuracy and computational scalability has produced methods that attempt to deal with this notorious 
trade-off.
The discussion that follows focuses specifically on this trade-off in the context of the \M model.
General approaches, such as those based on predictive processes \citep{banerjee2008gaussian} and those based on fixed-rank kriging \citep{cressie2008fixed}, will not be discussed; the interested reader is referred to the review of \citet{sun2012geostatistics}.


The computational complexity associated with the Mat{\'e}rn model is broadly governed by 
the input space dimension ($d$),  the number ($p$) of kernel parameters that must be estimated, and the number of data ($n$).
These challenges will be considered in turn.

First we consider the challenge of large $d$, 
which is often encountered in machine learning, when Gaussian process regression is performed on high-dimensional input spaces \citep{williams2006gaussian}. 
Since the Mat{\'e}rn model ${\cal M}_{\nu,\alpha}$ reproduces a Sobolev space 
up to an equivalent norm (c.f. Section \ref{SecNAaAT}), and $\nu>d/2$ is required to for elements of this space to be pointwise well-defined, it follows that $\nu$ must tend to infinity as $d$ tends to infinity, so that the \M model reduces to the Gaussian model \eqref{eqGauLim}.
The flexibility of some enhanced models is also lost in this limit; 
the condition $\mu \ge (d+1)/2+\kappa$ in the 
 $\widetilde{{\cal GW}}_{\kappa,\mu,\beta}$ model forces the parameter $\mu$ to go to infinity with $d$, which in turn forces $\widetilde{{\cal GW}}_{\kappa,\mu,\beta}$  to approach ${\cal M}_{\nu,\alpha}$.
 From this point of view the class  $\mathcal{GH}_{\kappa,\delta,\gamma,\beta}$ seems more promising to use for large $d$.
An additional remark is that, for $d \ge 5$, all Gaussian measures with Mat{\'e}rn covariance functions are orthogonal \citep{anderes2010}. This has philosophical consequences for Gaussian process regression when the \M model is viewed as a \emph{prior} distribution encoding \emph{a priori} belief, since a small change to the kernel parameters results in the entire support of the prior being changed.

Coupled to large input dimension $d$ is the challenge where there are a large number of parameters $p$ appearing in the model.
The multivariate Mat{\'e}rn model suffers from the fact, not only does $p$ increase exponentially with $d$, but the  conditions for validity of the model imply severe restrictions on the collocated correlation coefficient $\rho_{ij}$ in (\ref{stma}).
\citet{emery2022new} show that such restrictions become extremely severe already with $p=3$. 
Similar comments apply to other multivariate covariance functions, including the multivariate $\mathcal{GW}$ model in \citet{daley2015classes}.



Finally we consider the case where the number $n$ of data is large, entailing a $O(n^3)$ computational and $O(n^2)$ storage cost associated with the predictor \eqref{blup}.
Several approaches have been proposed to reduce these costs in the context of the \M model,
many of which take advantage of  the  (approximate) sparsity  of the covariance ($\bm{\Sigma}_n$) or precision ($\bm{\Sigma}_n^{-1})$, or its Cholesky factor ($\text{ch}(\bm{\Sigma}_n^{-1}))$:

\begin{itemize}
\item Sparsity in the covariance matrix $\bm{\Sigma}_n$ of the \M model  is directly exploited by enhanced versions of the \M model from Section \ref{surr}.
Such approaches can be useful when the (estimated) compact support is relatively small with respect to the  spatial extent of the sampling region, so that approximations are extremely sparse; see below for an empirical investigation of this point. 

\item The precision matrix $\bm{\Sigma}_n^{-1}$ associated with the Mat{\'e}rn model is in general non-sparse (except  for the case $d=1$ and $\nu=0.5$) but it turns out 
that the matrix values are in general relatively close to $0$,
i.e. $\bm{\Sigma}_n^{-1}$ is {\em quasi-sparse}.  As a consequence, approximating  
$\bm{\Sigma}_n^{-1}$ with a sparse matrix can be a good strategy.
A notable instance of this approach is the SPDE approach from Section \ref{SecSPDE}.
This  approach can be also motivated from  results in numerical linear algebra, which demonstrate that 
if the elements of a matrix show a property of decay, then the elements of its inverse also show a similar (and faster) behavior \citep{Benzi2016}.

\item Vecchia’s approximation \citep{vecchia1988} and its extensions
 \citep[e.g.][]{ datta2016hierarchical,guiness2018,ka2021,datta2021} 
imply a sparse approximation of of $\text{ch}(\bm{\Sigma}_n^{-1})$ and are often applied to the \M model, although they can be applied to any covariance model. One potential limitation of these method is that they depend on an ordering
of the variables and the choice of conditioning sets  which determines the Cholesky sparsity pattern
\citep[see][]{guiness2018}. 

\end{itemize}

\begin{table}[p]
\begin{center}
\Rotatebox{90}{\scalebox{0.775}{
\begin{tabular}{||c|c|c|c|c|c|c|c||c|c|c|c|c|c|c|c||c|c|c|c|c|c|c|c||}
\hline
\multicolumn{8}{||c||}{$\kappa=0$}                                                                                                                                                                                                                                                                                                                                                                                                                                                                 & \multicolumn{8}{c||}{$\kappa=1$}                                                                                                                                                                                                                                                                                                                                                                                                                                                                          & \multicolumn{8}{c||}{$\kappa=2$}                                                                                                                                                                                                                                                                                                                                                                                                                                                      \\ \hline
\multicolumn{1}{||c|}{$\mu$}                & \multicolumn{1}{c|}{$C$}                        & \multicolumn{2}{c|}{$\bm{\Sigma}_n$}                                                                                                    & \multicolumn{2}{c|}{$\bm{\Sigma}_n^{-1}$}                                                                                                & \multicolumn{2}{c||}{$\text{ch}(\bm{\Sigma}^{-1}_n)$}                                                                            & \multicolumn{1}{c|}{$\mu$}                         & \multicolumn{1}{c|}{$C$}                        & \multicolumn{2}{c|}{$\bm{\Sigma}_n$}                                                                                                    & \multicolumn{2}{c|}{$\bm{\Sigma}_n^{-1}$}                                                                                                & \multicolumn{2}{c||}{$\text{ch}(\bm{\Sigma}^{-1}_n)$}                                                                            & \multicolumn{1}{c|}{$\mu$}                & \multicolumn{1}{c|}{$C$}             & \multicolumn{2}{c|}{$\bm{\Sigma}_n$}                                                                                                    & \multicolumn{2}{c|}{$\bm{\Sigma}_n^{-1}$}                                                                                                & \multicolumn{2}{c||}{$\text{ch}(\bm{\Sigma}^{-1}_n)$}                                                                            \\ \hline
\multicolumn{1}{||c|}{}                     & \multicolumn{1}{c|}{}                                & \multicolumn{1}{c|}{\textit{$1156$}} & \multicolumn{1}{c|}{\textit{$4489$}} & \multicolumn{1}{c|}{\textit{$1156$}} & \multicolumn{1}{c|}{\textit{$4489$}} & \multicolumn{1}{c|}{\textit{$1156 $}} & \textit{$4489$} & \multicolumn{1}{c|}{}                              & \multicolumn{1}{c|}{}                                & \multicolumn{1}{c|}{\textit{$1156$}} & \multicolumn{1}{c|}{\textit{$4489$}} & \multicolumn{1}{c|}{\textit{$1156$}} & \multicolumn{1}{c|}{\textit{$4489$}} & \multicolumn{1}{c|}{\textit{$1156$}} & \textit{$4489$} & \multicolumn{1}{c|}{}                     & \multicolumn{1}{c|}{}                     & \multicolumn{1}{c|}{\textit{$1156$}} & \multicolumn{1}{c|}{\textit{$4489$}} & \multicolumn{1}{c|}{\textit{$1156$}} & \multicolumn{1}{c|}{\textit{$4489$}} & \multicolumn{1}{c|}{\textit{$1156$}} & \textit{$4489$} \\ \hline\hline

\multicolumn{1}{||c|}{$1.5$}                & \multicolumn{1}{c|}{$0.07$} & \multicolumn{1}{c|}{$98.4$}             & \multicolumn{1}{c|}{$98.3$}              & \multicolumn{1}{c|}{$32.3$}              & \multicolumn{1}{c|}{$ 1.46$}              & \multicolumn{1}{c|}{ $35.7$ }                   & \multicolumn{1}{c|}{ $2.17$ }                    & \multicolumn{1}{c|}{$2.5$}
&\multicolumn{1}{c|}{$ 0.11$}  
 & \multicolumn{1}{c|}{$97.1 $} & \multicolumn{1}{c|}{$96.7$}                                     & \multicolumn{1}{c|}{$3.80$}               & \multicolumn{1}{c|}{$1.67$}         & \multicolumn{1}{c|}{$ 6.39$}                                      & \multicolumn{1}{c|}{$4.44$}                                                                                  &
 
  \multicolumn{1}{c|}{$3.5$}                & \multicolumn{1}{c|}{$ 0.13$}              & \multicolumn{1}{c|}{$94.7 $}                                     & 
 \multicolumn{1}{c|}{$95.0$}                                      & \multicolumn{1}{c|}{$ 2.90$}                                      & \multicolumn{1}{c|}{$1.31 $}                                      & \multicolumn{1}{c|}{$5.63$}                                            &         \multicolumn{1}{c|}{$ 6.12$}                                         \\

\multicolumn{1}{||c|}{$4$}                  & \multicolumn{1}{c|}{$0.20$} & \multicolumn{1}{c|}{$90.1 $}             & \multicolumn{1}{c|}{$89.6 $}              & \multicolumn{1}{c|}{$ 45.9$}              & \multicolumn{1}{c|}{$56.0 $}              & \multicolumn{1}{c|}{$45.0$}          &  \multicolumn{1}{c|}{$54.9$}  & \multicolumn{1}{c|}{$4$}     
                      & \multicolumn{1}{c|}{$0.16 $}                         & \multicolumn{1}{c|}{$93 .4$}                                     & \multicolumn{1}{c|}{$ 93.3$}   & \multicolumn{1}{c|}{$34.2 $}                                      & \multicolumn{1}{c|}{$47.1 $}                                      & \multicolumn{1}{c|}{$36.6$}                                            &                   \multicolumn{1}{c|}{$48.7$}    
                                              & \multicolumn{1}{c|}{$4$}                  & \multicolumn{1}{c|}{$0.15$}              & \multicolumn{1}{c|}{$94.8$}                                     & \multicolumn{1}{c|}{$94.2$}                                      & \multicolumn{1}{c|}{$10.4 $}                                      & \multicolumn{1}{c|}{$ 11.9$}                                      & \multicolumn{1}{c|}{$15.0$}                                            &                          \multicolumn{1}{c|}{$22.6$}                        \\

\multicolumn{1}{||c|}{$8$}                  & \multicolumn{1}{c|}{$0.40$} & \multicolumn{1}{c|}{$66.4$}             & \multicolumn{1}{c|}{$65.7$}              & \multicolumn{1}{c|}{$56.6$}              & \multicolumn{1}{c|}{$64.2$}              & \multicolumn{1}{c|}{52.0}                    &  \multicolumn{1}{c|}{60.0}                    & \multicolumn{1}{c|}{$8$}                         & \multicolumn{1}{c|}{$0.28$}                         & \multicolumn{1}{c|}{$80.9$}                  & \multicolumn{1}{c|}{$80.7$}                      & \multicolumn{1}{c|}{$69.7$}                                      & \multicolumn{1}{c|}{$57.0$}                                      & \multicolumn{1}{c|}{$69.3$}                                      & \multicolumn{1}{c|}{$59.2$}                                                                                        &

\multicolumn{1}{c|}{$8$}                  & \multicolumn{1}{c|}{$0.25$}              & \multicolumn{1}{c|}{$84.7$}  & \multicolumn{1}{c|}{$84.4$}                           & \multicolumn{1}{c|}{$42.8$}                                      & \multicolumn{1}{c|}{$61.0$}                                      & \multicolumn{1}{c|}{$47.3$}                                            &        \multicolumn{1}{c|}{$64.8$}                                       \\

\multicolumn{1}{||c|}{$16$}                 & \multicolumn{1}{c|}{$0.80$}                         & \multicolumn{1}{c|}{$16.4 $}                                     & \multicolumn{1}{c|}{$ 15.2$}                                      & \multicolumn{1}{c|}{$59.8 $}                                      & \multicolumn{1}{c|}{$71.5 $}     & \multicolumn{1}{c|}{53.8}                    &          \multicolumn{1}{c|}{66.6}       

                             & \multicolumn{1}{c|}{$16$}                          & \multicolumn{1}{c|}{$0.54$}                         & \multicolumn{1}{c|}{$ 48.1$}                                     & \multicolumn{1}{c|}{$ 47.1$}                                      & \multicolumn{1}{c|}{$64.2 $}                                      & \multicolumn{1}{c|}{$80.1 $}                                      & \multicolumn{1}{c|}{64.0}                                            &                         \multicolumn{1}{c|}{76.0}                      & 
                             
                             \multicolumn{1}{c|}{$16$}                 & \multicolumn{1}{c|}{$0.45$}              & \multicolumn{1}{c|}{$ 59.1$}                                     & \multicolumn{1}{c|}{ $58.3 $}                                      & \multicolumn{1}{c|}{$ 51.7$}                                      & \multicolumn{1}{c|}{$ 73.0$}                                      & \multicolumn{1}{c|}{$54.6$}      &      \multicolumn{1}{c|}{$74.2$}                                             \\

\multicolumn{1}{||c|}{$32$}                 & \multicolumn{1}{c|}{$1.60$}                         & \multicolumn{1}{c|}{$0$}                                        & \multicolumn{1}{c|}{$0$}                                         & \multicolumn{1}{c|}{$ 61.2$}                                      & \multicolumn{1}{c|}{$79.0 $}                                      & \multicolumn{1}{c|}{55.8}          &      \multicolumn{1}{c|}{72.8}                                         &

\multicolumn{1}{c|}{$32$}                          & \multicolumn{1}{c|}{$1.04$}                         & \multicolumn{1}{c|}{$1.91$}                                        & \multicolumn{1}{c|}{$1.61$}                                         & \multicolumn{1}{c|}{$66.4$}                                      & \multicolumn{1}{c|}{$83.8 $}                                      & \multicolumn{1}{c|}{66.5}                                            &         \multicolumn{1}{c|}{82.2}  

                                  & \multicolumn{1}{c|}{$32$}                 & \multicolumn{1}{c|}{$0.86$}              & \multicolumn{1}{c|}{$11.0$}    & \multicolumn{1}{c|}{$10.1$}   &
        \multicolumn{1}{c|}{$52.3 $}                                      & \multicolumn{1}{c|}{$75.7 $}                                      & \multicolumn{1}{c|}{$54.8$}     & \multicolumn{1}{c|}{$77.2$}                                            \\

\multicolumn{1}{||c|}{$120$}                 & \multicolumn{1}{c|}{$6.01$}                         & \multicolumn{1}{c|}{$0$}                                        & \multicolumn{1}{c|}{$0$}                                         & \multicolumn{1}{c|}{$65.2 $}                                      & \multicolumn{1}{c|}{$ 76.6$}                                      & \multicolumn{1}{c|}{$58.9$}         &     \multicolumn{1}{c|}{$68.2$}     

                                   & \multicolumn{1}{c|}{$120$}                          & \multicolumn{1}{c|}{$3.82$}                         & \multicolumn{1}{c|}{$0$}        & \multicolumn{1}{c|}{$0$}                                    & \multicolumn{1}{c|}{$66.2$}                                         & \multicolumn{1}{c|}{$84.6 $}                                      & \multicolumn{1}{c|}{$ 65.9$}                                      & \multicolumn{1}{c|}{$81.2$}                                            &                                           \multicolumn{1}{c|}{$120$}                 &

                                    \multicolumn{1}{c|}{$3.09$}              & \multicolumn{1}{c|}{$0$}                                        & \multicolumn{1}{c|}{$0$}        & \multicolumn{1}{c|}{$52.4 $}                                      & \multicolumn{1}{c|}{$75.7 $}                                      & \multicolumn{1}{c|}{$54.1$}                                            &   \multicolumn{1}{c|}{$76.5$}                                            \\

\multicolumn{1}{||c|}{\textit{$\infty$}} & \multicolumn{1}{c|}{\textit{$\infty$}}            & \multicolumn{1}{c|}{$0$}                & \multicolumn{1}{c|}{$0$}                 & \multicolumn{1}{c|}{$66.2 $}                                      & \multicolumn{1}{c|}{$77.5$}                                      & \multicolumn{1}{c|}{$58.9$}               &                 \multicolumn{1}{c|}{$70.0$}                                   & \multicolumn{1}{c|}{\textit{$\infty$}}          & \multicolumn{1}{c|}{\textit{$\infty$}}            & \multicolumn{1}{c|}{$0$}                & \multicolumn{1}{c|}{$0$}                 & \multicolumn{1}{c|}{$66.2$}                                      & \multicolumn{1}{c|}{$84.9$}                                      & \multicolumn{1}{c|}{$65.1$}                                            &              \multicolumn{1}{c|}{$80.7$}                                  & \multicolumn{1}{c|}{\textit{$\infty$}} & \multicolumn{1}{c|}{\textit{$\infty$}} & \multicolumn{1}{c|}{$0$}                & \multicolumn{1}{c|}{$0$}                 & \multicolumn{1}{c|}{$52.4 $}                                      & \multicolumn{1}{c|}{$ 75.4$}                                      & \multicolumn{1}{c|}{$53.6$}                                            &         \multicolumn{1}{c|}{$74.6$}                                         \\ \hline
\end{tabular}
}}
 \end{center}
\caption{Sparsity (percentage of zero values in the upper triangular  part) of the covariance matrix $\bm{\Sigma}_n$, and quasi-sparsity (defined in the main text) in the precision matrix ($\bm{\Sigma}_n^{-1}$)
and its Cholesky 
factor ($\text{ch}(\bm{\Sigma}_n^{-1}$))
for the $\widetilde{{\cal GW}}_{\kappa,\mu,\beta}$  model.
The case
$\widetilde{{\cal GW}}_{\kappa,\infty,\beta}$
 corresponds to the Mat{\'e}rn  model   ${\cal M}_{\nu+1/2,\beta}$. The $\beta$ parameters are chosen so that
the practical range of the  Mat{\'e}rn  model is equal to $0.15$.
} 
\label{tab2222}
\end{table}

\begin{table}[htbp] 
\begin{center}
\Rotatebox{90}{\scalebox{0.8}{
\begin{tabular}{||c|c|c|c|c|c|c|c||c|c|c|c|c|c|c|c||c|c|c|c|c|c|c|c||}
\hline
\multicolumn{8}{||c||}{$\kappa=0$}                                                                                                                                                                                                                                                                                                                                                                                                                                                                 & \multicolumn{8}{c||}{$\kappa=1$}                                                                                                                                                                                                                                                                                                                                                                                                                                                                          & \multicolumn{8}{c||}{$\kappa=2$}                                                                                                                                                                                                                                                                                                                                                                                                                                                      \\ \hline
\multicolumn{1}{||c|}{$\mu$}                & \multicolumn{1}{c|}{$C$}                        & \multicolumn{2}{c|}{$\bm{\Sigma}_n$}                                                                                                    & \multicolumn{2}{c|}{$\bm{\Sigma}_n^{-1}$}                                                                                                & \multicolumn{2}{c||}{$\text{ch}(\bm{\Sigma}^{-1}_n)$}                                                                            & \multicolumn{1}{c|}{$\mu$}                         & \multicolumn{1}{c|}{$C$}                        & \multicolumn{2}{c|}{$\bm{\Sigma}_n$}                                                                                                    & \multicolumn{2}{c|}{$\bm{\Sigma}_n^{-1}$}                                                                                                & \multicolumn{2}{c||}{$\text{ch}(\bm{\Sigma}^{-1}_n)$}                                                                            & \multicolumn{1}{c|}{$\mu$}                & \multicolumn{1}{c|}{$C$}             & \multicolumn{2}{c|}{$\bm{\Sigma}_n$}                                                                                                    & \multicolumn{2}{c|}{$\bm{\Sigma}_n^{-1}$}                                                                                                & \multicolumn{2}{c||}{$\text{ch}(\bm{\Sigma}^{-1}_n)$}                                                                            \\ \hline
\multicolumn{1}{||c|}{}                     & \multicolumn{1}{c|}{}                                & \multicolumn{1}{c|}{\textit{$1156$}} & \multicolumn{1}{c|}{\textit{$4489$}} & \multicolumn{1}{c|}{\textit{$1156$}} & \multicolumn{1}{c|}{\textit{$4489$}} & \multicolumn{1}{c|}{\textit{$1156 $}} & \textit{$4489$} & \multicolumn{1}{c|}{}                              & \multicolumn{1}{c|}{}                                & \multicolumn{1}{c|}{\textit{$1156$}} & \multicolumn{1}{c|}{\textit{$4489$}} & \multicolumn{1}{c|}{\textit{$1156$}} & \multicolumn{1}{c|}{\textit{$4489$}} & \multicolumn{1}{c|}{\textit{$1156$}} & \textit{$4489$} & \multicolumn{1}{c|}{}                     & \multicolumn{1}{c|}{}                     & \multicolumn{1}{c|}{\textit{$1156$}} & \multicolumn{1}{c|}{\textit{$4489$}} & \multicolumn{1}{c|}{\textit{$1156$}} & \multicolumn{1}{c|}{\textit{$4489$}} & \multicolumn{1}{c|}{\textit{$1156$}} & \textit{$4489$} \\ \hline\hline

\multicolumn{1}{||c|}{$1.5$}                & \multicolumn{1}{c|}{$0.20$} & \multicolumn{1}{c|}{$90.0$}             & \multicolumn{1}{c|}{$89.1$}              & \multicolumn{1}{c|}{$0$}              & \multicolumn{1}{c|}{$0$}              & \multicolumn{1}{c|}{ $0$ }                   & \multicolumn{1}{c|}{ $0$ }                  

  & \multicolumn{1}{c|}{$2.5$}
&\multicolumn{1}{c|}{$0.28$}  
 & \multicolumn{1}{c|}{$80.5 $} & \multicolumn{1}{c|}{$80.2$}                                     & \multicolumn{1}{c|}{$0$}               & \multicolumn{1}{c|}{$0$}         & \multicolumn{1}{c|}{$0 $}                                      & \multicolumn{1}{c|}{$0$}                                                                                  &
 
  \multicolumn{1}{c|}{$3.5$}                & \multicolumn{1}{c|}{$ 0.35$}              & \multicolumn{1}{c|}{$ 72.6$}                                     & 
 \multicolumn{1}{c|}{$71.4$}                                      & \multicolumn{1}{c|}{$ 0$}                                      & \multicolumn{1}{c|}{$ 0$}                                      & \multicolumn{1}{c|}{$0$}                                            &         \multicolumn{1}{c|}{$0 $}                                         \\

\multicolumn{1}{||c|}{$4$}                  & \multicolumn{1}{c|}{$0.53$} & \multicolumn{1}{c|}{$48,7 $}             & \multicolumn{1}{c|}{$47.3$}              & \multicolumn{1}{c|}{$ 1.03$}              & \multicolumn{1}{c|}{$ 9.62$}              & \multicolumn{1}{c|}{$0.58$}          &  \multicolumn{1}{c|}{$2.77$}  &

 \multicolumn{1}{c|}{$4$}     
                      & \multicolumn{1}{c|}{$0.42 $}                         & \multicolumn{1}{c|}{$65.0$}                                     & \multicolumn{1}{c|}{$  63.7$}   & \multicolumn{1}{c|}{$0 $}                                      & \multicolumn{1}{c|}{$ 1.03$}                                      & \multicolumn{1}{c|}{$1.12$}                                            &                   \multicolumn{1}{c|}{$0.81$}    
                                              & \multicolumn{1}{c|}{$4$}                  & \multicolumn{1}{c|}{$0.39$}              & \multicolumn{1}{c|}{$67.2$}                                     & \multicolumn{1}{c|}{$66.4$}                                      & \multicolumn{1}{c|}{$ 0$}                                      & \multicolumn{1}{c|}{$0$}                                      & \multicolumn{1}{c|}{$0$}                                            &                          \multicolumn{1}{c|}{$0$}                        \\

\multicolumn{1}{||c|}{$8$}                  & \multicolumn{1}{c|}{$1.07$} & \multicolumn{1}{c|}{$1.43$}             & \multicolumn{1}{c|}{$1.21$}              & \multicolumn{1}{c|}{$4.14$}              & \multicolumn{1}{c|}{$25.6$}              & \multicolumn{1}{c|}{$1.84$}                    &  \multicolumn{1}{c|}{$8.15$}                    
& \multicolumn{1}{c|}{$8$}                         & \multicolumn{1}{c|}{$0.75$}                         & \multicolumn{1}{c|}{$20.7$}                  & \multicolumn{1}{c|}{$19.7$}                      & \multicolumn{1}{c|}{$5.91$}                                      & \multicolumn{1}{c|}{$14.5$}                                      & \multicolumn{1}{c|}{$4.49$}                                      & \multicolumn{1}{c|}{$12.9$}                                                                                        & 
\multicolumn{1}{c|}{$8$}                  & \multicolumn{1}{c|}{$0.67$}              & \multicolumn{1}{c|}{$30.7$}  & \multicolumn{1}{c|}{$29.6$}                           & \multicolumn{1}{c|}{$0.71$}                                      & \multicolumn{1}{c|}{$2.73$}                                      & \multicolumn{1}{c|}{$2.13$}                                            &        \multicolumn{1}{c|}{$6.90$}                                       \\

\multicolumn{1}{||c|}{$16$}                 & \multicolumn{1}{c|}{$2.13$}                         & \multicolumn{1}{c|}{$0$}                                     & \multicolumn{1}{c|}{$0 $}                                      & \multicolumn{1}{c|}{$ 15.0$}                                      & \multicolumn{1}{c|}{$43.5$}     & \multicolumn{1}{c|}{$4.98$}                    &          \multicolumn{1}{c|}{$21.3$}       

                             & \multicolumn{1}{c|}{$16$}                          & \multicolumn{1}{c|}{$1.43$}                         & \multicolumn{1}{c|}{$ 0$}                                     & \multicolumn{1}{c|}{$ 0$}                                      & \multicolumn{1}{c|}{$ 21.2$}                                      & \multicolumn{1}{c|}{$ 45.2$}                                      & \multicolumn{1}{c|}{$17.1$}                                            &                         \multicolumn{1}{c|}{$30.8$}                      & 
                             
                             \multicolumn{1}{c|}{$16$}                 & \multicolumn{1}{c|}{$1.21$}              & \multicolumn{1}{c|}{$1.72$}                                     & \multicolumn{1}{c|}{ $ 1.29$}                                      & \multicolumn{1}{c|}{$ 13.7$}                                      & \multicolumn{1}{c|}{$ 25.8$}                                      & \multicolumn{1}{c|}{$14.8$}      &      \multicolumn{1}{c|}{$27.9$}                                             \\

\multicolumn{1}{||c|}{$32$}                 & \multicolumn{1}{c|}{$4.27$}                         & \multicolumn{1}{c|}{$0$}                                        & \multicolumn{1}{c|}{$0$}                                         & \multicolumn{1}{c|}{$ 22.8$}                                      & \multicolumn{1}{c|}{$44.4$}                                      & \multicolumn{1}{c|}{$12.7$}          &      \multicolumn{1}{c|}{$24.3$}                                         &

\multicolumn{1}{c|}{$32$}                          & \multicolumn{1}{c|}{$2.78$}                         & \multicolumn{1}{c|}{$0$}                                        & \multicolumn{1}{c|}{$0$}                                         & \multicolumn{1}{c|}{$37.2$}                                      & \multicolumn{1}{c|}{$55.0 $}                                      & \multicolumn{1}{c|}{$29.0$}                                            &         \multicolumn{1}{c|}{$41.0$}  

                                  & \multicolumn{1}{c|}{$32$}                 & \multicolumn{1}{c|}{$2.29$}              & \multicolumn{1}{c|}{$0$}    & \multicolumn{1}{c|}{$0$}   &
        \multicolumn{1}{c|}{$25.5$}                                      & \multicolumn{1}{c|}{$23.5$}                                      & \multicolumn{1}{c|}{$24.8$}     & \multicolumn{1}{c|}{$42.7$}                                            \\

\multicolumn{1}{||c|}{$120$}                 & \multicolumn{1}{c|}{$16.0$}                         & \multicolumn{1}{c|}{$0$}                                        & \multicolumn{1}{c|}{$0$}                                         & \multicolumn{1}{c|}{$21.0$}                                      & \multicolumn{1}{c|}{$ 46.2$}                                      & \multicolumn{1}{c|}{$8.33$}         &     \multicolumn{1}{c|}{$21.9$}     

                                   & \multicolumn{1}{c|}{$120$}                          & \multicolumn{1}{c|}{$10.2$}                         & \multicolumn{1}{c|}{$0$}        & \multicolumn{1}{c|}{$0$}                                    & \multicolumn{1}{c|}{$40.8$}                                         & \multicolumn{1}{c|}{$59.8 $}                                      & \multicolumn{1}{c|}{$ 32.2$}                                      & \multicolumn{1}{c|}{$46.2$}                                            &                                           \multicolumn{1}{c|}{$120$}                 &
                                    \multicolumn{1}{c|}{$8.24$}              & \multicolumn{1}{c|}{$0$}                                        & \multicolumn{1}{c|}{$0$}        & \multicolumn{1}{c|}{$26.0$}                                      & \multicolumn{1}{c|}{$10.9 $}                                      & \multicolumn{1}{c|}{$25.9$}                                            &   \multicolumn{1}{c|}{$46.4$}                                            \\

\multicolumn{1}{||c|}{\textit{$\infty$}} & \multicolumn{1}{c|}{\textit{$\infty$}}            & \multicolumn{1}{c|}{$0$}                & \multicolumn{1}{c|}{$0$}                 & \multicolumn{1}{c|}{$23.4$}                                      & \multicolumn{1}{c|}{$47.9$}                                      & \multicolumn{1}{c|}{$9.85$}               &                 \multicolumn{1}{c|}{$24.0$}                                   & \multicolumn{1}{c|}{\textit{$\infty$}}          & \multicolumn{1}{c|}{\textit{$\infty$}}            & \multicolumn{1}{c|}{$0$}                & \multicolumn{1}{c|}{$0$}                 & \multicolumn{1}{c|}{$42.6$}                                      & \multicolumn{1}{c|}{$61.3$}                                      & \multicolumn{1}{c|}{$34.3$}                                            &              \multicolumn{1}{c|}{$48.0$}                                  & \multicolumn{1}{c|}{\textit{$\infty$}} & \multicolumn{1}{c|}{\textit{$\infty$}} & \multicolumn{1}{c|}{$0$}                & \multicolumn{1}{c|}{$0$}                 & \multicolumn{1}{c|}{$26.0$}                                      & \multicolumn{1}{c|}{$ 18.2$}                                      & \multicolumn{1}{c|}{$25.9$}                                            &         \multicolumn{1}{c|}{$46.6$}                                         \\ \hline
\end{tabular}
}}
 \end{center}
\caption{As in Table \ref{tab2222}, but with $\beta$ chosen such that the practical range of the \M model is equal to 0.4.
 } \label{tab333}
\end{table}

\noindent It is instructive to numerically investigate the sparseness of matrices associated with enhancements of the \M model, and for this we focus on the $\widetilde{{\cal GW}}_{\kappa,\mu,\beta}$ model,
 which allows us to  switch from a model with compact support of radius 
 $$
 C=
\beta \left(\frac{\Gamma(\mu+2\kappa +1)} {\Gamma(\mu)}\right)^{\frac{1}{1+2\kappa}}
 $$
 to the Mat{\'e}rn model by increasing the $\mu$  parameter.
In our experiment, the sparseness of $\bm{\Sigma}_n$   and 
the \emph{quasi-sparseness} of $\bm{\Sigma}_n^{-1}$ and $\text{ch}(\bm{\Sigma}_n^{-1})$ 
are reported, the latter being defined as the percentage of values in the upper triangular matrix with  absolute value lower than an arbitrary
small constant $\epsilon$, and  in our example we set $\epsilon=1.e-8$.

The empirical assessment considers $n=1,156$ and $n=4,489$ location sites over $[0,1]^2$, where the  points are equally spaced by $0.03$ and $0.015$ respectively in a regular grid.
For $\nu=0,1,2$, we set $\beta$ such that  the  practical range of the Mat{\'e}rn model  is equal to $0.15$ ($\beta=0.050, 0.0316, 0.0253$ respectively),
and consider increasing $\mu=1.5+\kappa,4, 8, 16, 32, 120,  \infty$
(with $\widetilde{{\cal GW}}_{\kappa,\infty,\beta}$ being the Mat{\'e}rn model
  ${\cal M}_{\kappa+1/2,\beta})$. 

  The results are reported in Table  \ref{tab2222}. 
  For the low values $\mu=1.5, 2.5, 3.5$ and $\nu=0, 1, 2$,
  the covariance matrix is highly sparse, while the sparseness decreases when increasing $\mu$, as expected.
  There is a clear
trade-off  between the sparseness of  $\bm{\Sigma}_n$ and quasi-sparseness of $\bm{\Sigma}_n^{-1}$ and $\text{ch}(\bm{\Sigma}_n^{-1})$
  for each $\nu=0, 1, 2$. 
  However, when increasing $\mu$,  that is when $\bm{\Sigma}_n$ approaches the  Mat{\'e}rn covariance matrix, then
  $\bm{\Sigma}_n^{-1}$ or  $\text{ch}(\bm{\Sigma}_n^{-1})$ tends to be  highly quasi-sparse.
  

We replicate the same experiment but with a practical range of the Mat{\'e}rn  model equal to $0.4$. This leads to $\beta=0.133, 0.084, 0.067$ for $\nu=0, 1, 2$ respectively.
  The results are reported in Table  \ref{tab333}.
The conclusions are the same of the previous setting but in this case, we have lower levels of sparseness   for $\bm{\Sigma}_n$ and of quasi-sparseness  for $\bm{\Sigma}_n^{-1}$
and $\text{ch}(\bm{\Sigma}_n^{-1}$.

These numerical experiments highlight a clear trade-off between the (quasi-)sparseness of
 $\bm{\Sigma}^{-1}_n$ (or $\text{ch}(\bm{\Sigma}_n^{-1})$)  and  $\bm{\Sigma}_n$ 
 when increasing $\mu$ for fixed $\beta$ and $\nu$
 i.e. when switching from a compactly supported to a globally supported Mat{\'e}rn model.
In particular, when $\mu \to \infty$ (the  Mat{\'e}rn model), then  $\bm{\Sigma}_n^{-1}$  
is highly quasi-sparse and   $\bm{\Sigma}_n$ is dense.
In contrast, 
when $\mu$ is small then  $\bm{\Sigma}_n^{-1}$  
is not quasi-sparse yet $\bm{\Sigma}_n$ is highly sparse.
 This seems to suggest that sparse precision matrix approximation should work reasonably well for the Mat{\'e}rn model, 
but could be problematic when handling  data  exhibiting short compactly supported dependence.
In this case a better approach should be to exploit the sparsity of $\bm{\Sigma}_n$, as enabled by enhanced versions of the \M model.

\section{Conclusion}

The impact of the \M model since its conception has been substantial, and the model continues to be widely used, across a broad range of scientific disciplines and beyond.
While the original motivation for the Mat{\'e}rn model came from its flexibility in  context of spatial interpolation, there is now also a rich literature of alternative and enhanced versions of the model. 
In particular, the SPDE and related approaches enable one to define analogues of the \M model on quite general domains, admitting sparse approximations to precision matrices, while recent advances in enhanced models with compact support can facilitate scalable computation through sparse approximation of covariance matrices, and are well-suited to processes with short-scale dependence.
The theoretical and empirical properties of these enhanced models have been recently and actively studied.
On the other hand, there remain open theoretical issues of practical importance, such as parameter estimation at finite sample sizes, and the impact of parameter estimation on the performance of the associated predictions.

Our current understanding of the \M model has emerged as the result of engagement between scientists and practitioners from different disciplines, and our hope is that this multi-disciplinarity perspective will shine further light onto the \M model. 





\bibliographystyle{apalike}


\appendix 
\onecolumn

\begin{center}
\textbf{\Large Supplementary Material}
\end{center}

\section{Modelling through \M 
in unconventional scenarios: the extended version} \label{Appendix}

One might object that the Mat{\'e}rn class is limited to scalar-valued random fields that are stationary and isotropic. While this being true, it is also true that the Mat{\'e}rn class represents the building block for way more sophisticated scenarios. We list them throughout. 
\begin{enumerate}
\item Scalar-valued random fields.
\begin{enumerate}
\item {\em Anisotropies. } If spatial dimension develops over preferential directions, isotropy is no longer a realistic assumption for spatial modeling. Several types of anisotropies are challenged in \cite{allard2016anisotropy}, and it is shown that the Mat{\'e}rn class can be composed with {\em ad hoc} deformations so to take into account preferential directions in terms of spatial dependence.

\item {\em Nonstationarity.} The Mat{\'e}rn kernel has been used by \cite{paciorek2006spatial} to build nonstationary models. 
Consider a collection of Gaussian distributions indexed by their mean, such that the element with mean $\bm{x}$ has covariance matrix $\Sigma_{\bm{x}}$. 
Let $$Q_{\bm{x},\bm{y}}= (\bm{x}-\bm{y})^{\top} \left ( \frac{\Sigma_{\bm{x}}+\Sigma_{\bm{y}}}{2}\right )^{-1}(\bm{x}-\bm{y}) $$
and 
\begin{eqnarray*}
K(\bm{x},\bm{y})&=& \Big | \Sigma_{\bm{x}}\Big |^{1/4} \Big | \Sigma_{\bm{y}}\Big |^{1/4} \Bigg | \frac{\Sigma_{\bm{x}}+\Sigma_{\bm{y}}}{2}\Bigg |^{-1/2}  {\cal M}_{\nu,\alpha} \left ( \sqrt{Q_{\bm{x},\bm{y}}} \right ). 
\end{eqnarray*}
Then, $K$ is positive definite. 
This approach was recently generalised in \citet{roininen2019hyperpriors}.
Alternatively, one may induce non-stationarity by \emph{warping} the inputs of the Mat\'{e}rn covariance function,  
as $K(\bm{x},\bm{y}) = \mathcal{M}_{\nu,\alpha}(\|\bm{w}(\bm{x})-\bm{w}(\bm{y})\|)$ for some diffeomorphism $\bm{w} : \R^d \rightarrow \R^d$ that may itself be parametrised.
The case in which $\bm{w}$ is pa\-ra\-me\-trised by a deep neural network was explored in \citet{wilson2016deep}.

\item {\em Graphs and Quasi Metric spaces.} \cite{anderes2020isotropic} consider graphs with Euclidean edges, equipped with either the geodesic distance over the graph, or the resistance metric. They prove that ${\cal M}_{\nu,\alpha}$ can be composed with the resistance metric over the graph provided $0<\nu\le 1/2$. More recently, \cite{menegatto2020gneiting} provide a generalisation of this setting by considering quasi-metric spaces. Apparently, similar restrictions hold for this case. Recently, \cite{https://doi.org/10.48550/arxiv.2205.06163} adopt a different approach to build random fields with their covariance structure on metric graphs. 
Space-time version of the Mat{\'e}rn class, for space being a graph with Euclidean edges, have been considered by \cite{tang2020space} and by \cite{porcu2022nonseparable}.

\item {\em Space-time.} For a space-time Gaussian random field $\{Z(\bm{x},t), \; \bm{x} \in \R^{d+1}, t \in \R \}$, a typical second-order assumption is that the covariance is isotropic in space and stationary over time. That is, 
\begin{equation}
\label{space-time}
\begin{array}{rcl} 
 {\rm Cov} \left ( Z(\bm{x},t),Z(\bm{y},t')\right ) =K \left ( \|\bm{x}-\bm{y}\|,|t-t'|\right )
 \end{array}
\end{equation}
 for all $(\bm{x},t),(\bm{y},t') \in \R^d \times \R$. 
 The Mat{\'e}rn function has been used as a building block for such a structure. 
Of particular interest are \emph{non-separable} covariance functions, that allow for an interaction between space and time, and in this context \cite{gneiting1} and \cite{zast-porcu} prove that 
\begin{equation}\label{maast}
K(x,u)= \frac{\sigma^2}{\psi(u^2)} {\cal M}_{\nu,\alpha} \left ( \frac{x}{\psi(u^2)}\right ), \;x,u \ge 0, 
\end{equation}
generate a valid space-time covariance function of the type (\ref{space-time}). Here, $\psi$ is a strictly positive function having a completely monotonic derivative \citep{berg2008stieltjes}.

\item {\em Non Gaussian Fields.}  
Through consideration of transformations $w : \R \rightarrow \R$, one can use the Mat\'{e}rn model as the basis for a range of non-Gaussian models $\tilde{Z}(\bm{x}) = w(Z(\bm{x}))$.
However, the covariance function of $\tilde{Z}$ will not be Mat\'{e}rn in general.
The question of whether there exist non-Gaussian processes whose covariance function is nevertheless of Mat\'{e}rn class was answered positively in \cite{aaberg2011class}. 
 \cite{yan2018gaussian} have proposed trans-Gaussian random fields with Mat{\'e}rn covariance function. \cite{bolin2014spatial} and subsequently \cite{wallin2015geostatistical} have provided SPDE based constructions for non Gaussian Mat{\'e}rn fields.  A general class of non Gaussian fields  with kernel $g({\cal M}_{\nu,\alpha})$ with $g(\cdot)$ a suitable function that preserves positive definiteness can be obtained through  a transformation of (independent replicates of)  a Gaussian field (see for instance \cite{Palacios:Steel:2006,Xua:Genton:2017,Bevilacqua_et_al:2021,mmaa}).

\item {\em Mixtures of Gaussian Fields.}  
One can model a real-world phenomenon as a superposition $Z_1 + Z_2$ where $Z_1$ might be selected to capture an overall trend (e.g. $Z_1(\bm{x}) = c_1 \phi_1(\bm{x}) + \dots + c_p \phi_p(\bm{x})$ for some fixed functions $\phi_1,\dots,\phi_p$ and some coefficients $\bm{c} \in \R^p$ and $Z_2$ could be Mat\'{e}rn, to capture the level of smoothness of the process.
Such processes appear in Bayesian linear regression, where the coefficient $\bm{c}$ is viewed as random due to epistemic uncertainty, for example $\bm{c} \sim \mathcal{N}(\bm{\mu},\bm{\Sigma})$, and both $\bm{c}$ and $Z_2$ are to be jointly inferred.
The resulting mixture $Z_1 + Z_2$ is then again a Gaussian random field.
\item {\em Mat{\'e}rn Models with Modified Tails.}
The Mat\'ern covariance function decays exponentially with distance. This can be a drawback in the presence of long memory. Some classes of covariance functions allow to index long versus short memory in spatial data. The generalised Cauchy \citep{gneiting2004stochastic} and the Dagum \citep{berg2008dagum} are prominent examples of families allowing to index both fractal dimensions and long memory, also termed Hurst effect. Yet, these families do not allow to parameterise the smoothness in the same fashion of the Mat{\'e}rn class. This dilemma has preoccupied several scientists. Below we describe the approaches devoted to modify the tails 
of the Mat{\'e}rn model while preserving (a) positive definiteness and (b) local behavior, in turn connected with mean square differentiability and Sobolev space parameterisation. Each of the contributions below has different motivations as explained throughout. 
\begin{enumerate}
\item[(a)] {\em Mat\'ern models with periodic tails.} A modification of the spectral density of the Mat\'ern model has been proposed by \cite{laga2017modified}.   The primary idea is to propose an isotropic class of spectral densities, $\widetilde{M}_{\nu,\alpha,\sigma^2}$, that is connected to the Mat{\'e}rn family $\widehat{M}_{\nu,\alpha,\sigma^2}$ through the identity
\begin{equation}
\label{laga}
   \widetilde{M}_{\nu,\alpha,\sigma^2}(z)=\left( b^2+ z^2 \right)^\xi  \widehat{\mathcal{M}}_{\nu,\alpha,\sigma^2}(z),
\end{equation}
with $b\geq 0$ and $\xi < \nu$. While $b$ is an additional range parameter, $\xi$ is related to the smoothness of the respective process. More precisely, the random field is $k$-times mean square differentiable if and only if $\nu-\xi > k$. In the limit case $\xi\rightarrow 0$, the traditional Mat\'ern model is recovered. 

This scale in the spectral density produces a shift in the mode of the spectrum; thus, it is particularly useful to obtain processes with strong periodicities. The covariance function associated to (\ref{laga}) does not have a known explicit expression, so statistical methodologies in the spectral domain should be employed when dealing with this model.

\item[(b)] {\em Hybrid Models.}
Another generalisation of the Ma\-t\'ern model exploits the fact that it can be written as a scale mixture of a Gaussian kernel against a probability density function  
\begin{equation}
    \label{mixture}
    \mathcal{M}_{\nu,\alpha}(x) = \int_0^\infty \exp(-u x^2) \pi_{IG}(u;\nu,1/(4\alpha^2) ) \text{d}u
\end{equation}
of the inverse gamma type. 
\cite{alegria2023hybrid} make use of the identity above to create {\em hybrid} models that allow to preserve the local properties of the conventional Mat{\'e}rn model while attaining more flexible behaviours at large distances. 
 One potential hybrid construction, called Mat{\'e}rn-Cauchy and denoted $\mathcal{MC}_{\nu_1,\nu_2,\alpha,\xi}$, is attained through
\begin{equation}
    \label{mc}    \mathcal{MC}_{\nu_1,\nu_2,\alpha,\xi}(x)
    =  \int_0^\xi  \exp(-u x^2) \pi_G(u;\nu_1/2,\alpha) \text{d}u  +  \int_\xi^\infty \exp(-u x^2) \pi_{IG}(u;\nu_2,1/(4\alpha^2) ) \text{d}u.
\end{equation}
Here, $\pi_G$ is the gamma probability density function,  with a shape-rate parameterisation, $\pi_{IG}$ is the density in (\ref{inverse_gamma}), $\alpha > 0$ is a range parameter, $\nu_1 > 0$ controls the polynomial rate of decay of the covariance, $\nu_2 > 0$ indexes the mean square differentiability, and $\xi \geq 0$ is an additional parameter that balances the Mat\'ern and Cauchy contributions to the total covariance function. As $\xi\rightarrow 0$, the hybrid model tends to a Mat\'ern covariance.  
A closed form expression for this model is provided in  \cite{alegria2023hybrid}. Note that this covariance function is positive definite in any dimension, and is a natural competitor of the model (\ref{anciano}).

Another hybrid model is constructed in \cite{alegria2023hybrid} by replacing the Gaussian kernel in (\ref{mixture}) with a difference of Gaussian kernels 
$$\left(a \exp(-u b x^2) -  \exp(-u x^2)\right) (a-1)^{-1},$$ where $a$ and $b$ satisfy the condition $1<b < a^{2/d}$, in order to obtain a positive definite kernel in $\mathbb{R}^d$ \citep{posa2022special}. The resulting model has a local behavior of Mat\'ern type and allows for negative correlations at large distances (hole effect). The parameters $a$ and $b$ control the sharpness of the hole effect. When $a$ is arbitrarily large, the conventional Mat\'ern model is recovered. Algebraically closed expressions are reported in \cite{alegria2023hybrid}.

\end{enumerate}

\end{enumerate}

\item Vector-valued random fields. \\
Let $\{\bm{Z}(\bm{x}), \; \bm{x} \in \R^d\} \subset \R^p$ be a $p$-variate random field with isotropic covariance mapping $\bm{K}: \R^d \to \R^{p \times p}$ having elements $K_{ij}$ defined as 
$$K_{ij}(\bm{x})= {\rm cov} \left ( Z_i(\bm{0}),Z_j(\bm{x})\right ), \; \bm{0}, \bm{x} \in \R^d. $$
Specifically, $K_{ii}$ and $K_{ij}$ are termed auto and cross-covariance function.

\begin{enumerate}
\item {\em Multivariate Spatial.} There has been a ple\-tho\-ra of approaches related to multivariate spatial modeling, and the reader is referred to \cite{Genton:Kleiber:2014}. \cite{Gneiting:Kleibler:Schlather:2010} have proposed a multivariate covariance structure of the type
\begin{equation*}
K_{ij}(\bm{x}) = \sigma_{ii} \sigma_{jj} \rho_{ij} {\cal M}_{\nu_{ij},\alpha_{ij}} (\|\bm{x}\|), \; \bm{x} \in \R^d,
\end{equation*}
where $\sigma_{ii}^2$ is the variance of $Z_i$ from $\bm{Z}$, and $\rho_{ij}$ is the collocated correlation coefficient. There are restrictions on the parameters $\nu_{ij}, \alpha_{ij}$ and $\rho_{ij}$ to preserve positive definiteness, and often the restrictions on the collocated correlations coefficients $\rho_{ij}$ are severe. This motivated alternative approaches to alleviate the parametric restrictions, and the reader is referred to \cite{Apanasovich} and more recently to \cite{emery2022new}.

\item {\em Multivariate space-time.} The setting above can be generalised to space-time by considering $\{ \bm{Z}(\bm{x},t), \; \bm{x} \in \R^d, t\in \R  \} \subset \R^p$ a $p$-variate random field with isotropic covariance mapping $\bm{K}: \R^d \times \R \to \R^{p \times p}$ having elements $K_{ij}$ defined as 
$$K_{ij}(\bm{x},u)= {\rm cov} \left ( Z_i(\bm{0},t),Z_j(\bm{x},t+u)\right ) 
$$
for $\bm{0}, \bm{x} \in \R^d$, $t,u \in \R$. \cite{bourotte2016flexible} consider mappings $K_{ij}$ of the type
$$ K_{ij}(\bm{x},u) = \frac{\sigma_{ii} \sigma_{jj} \rho_{ij}}{\psi(u^2)^{d/2}} {\cal M}_{\nu_{ij},\alpha_{ij} \psi(u^2)} (\|\bm{x}\|)
$$
for $\bm{x} \in \R^d$, $u \in \R$ and  
a suitable positive valued and continuous function $\psi$. This setting has been recently generalised
by \cite{allard2022fully} and through a technical approach by \cite{porcu2022criteria}: for both contribution, the idea is to replace (pointwise) the mapping $\psi$ with the mapping $\bm{\psi}$ having continuous and strictly 
positive elements $\psi_{ij}$. 

\item {\em Multivariate Nonstationary.} \cite{kleiber2012nonstationary} derive a class of matrix valued covariance functions where the direct and cross-covariance functions belong to the Matérn class. The parameters of the Matérn class are allowed to vary with location, yielding local variances, local ranges, local geometric anisotropies and local smoothnesses. 

Define $\Sigma_{i,\bm{x}} = \Sigma_i(\bm{x}): \R^
d \to \R^{d \times d}$ and assume 
$\Sigma_i(\bm{x})$ is positive definite for all $i$ and all $\bm{x} \in \R^d$.
Let 
$$Q_{ij;\bm{x},\bm{y}}= (\bm{x}-\bm{y})^{\top} \left ( \frac{\Sigma_{i,\bm{x}}+\Sigma_{j,\bm{y}}}{2}\right )^{-1}(\bm{x}-\bm{y})
$$
and
$$
K_{ij}(\bm{x},\bm{y})= \rho_{ij}\sigma_{i,\bm{x}} \sigma_{j,\bm{y}} {\cal M}_{\nu,\alpha} \left ( \sqrt{Q_{ij;\bm{x},\bm{y}}} \right ). $$
Then, $\bm{K}(\cdot,\cdot) = [K_{ij}(\cdot,\cdot)]_{i,j=1}^p$ is positive definite.

\item {\em Multivariate Mat{\'e}rn with Dimple.}  In a bivariate spatial context, each element of the matrix-valued Mat\'ern covariance function admits a scale-mixture representation as in (\ref{mixture}). \cite{alegria2021bivariate} considered a modification of such a mixture to obtain a generalization of the Mat\'ern model, given by  
\begin{eqnarray}
    \label{dimple}
    \widetilde{K}_{ij}(x) =   \sigma_{ii}\sigma_{jj}\rho_{ij}  \int_0^\infty \exp(-u x^2) 
  g_{ij}(u;\xi) \times \pi_{IG}(u;\nu_{ij},1/(4\alpha_{ij}^2) ) \text{d}u,
\end{eqnarray}
where $g_{ii}(u;\xi) = 1$, and $g_{ij}(u;\xi) = 1\{u \leq \xi\} -  1\{u \geq \xi\}$ for $i\neq j$, with $1\{\cdot\}$ being the indicator function and $\xi$ a nonnegative parameter, and where $\pi_{IG}$ is the probability density function of an inverse gamma random variable, that is
\begin{equation}
    \label{inverse_gamma}
    \pi_{IG}(z; a,b) = \frac{b^a}{\Gamma(a)}z^{-a-1} \exp(-b/z), \;z>0. 
\end{equation}

Observe that the diagonal elements of the matrix-valued covariance are not altered; thus, the  appealing local attributes of the Mat\'ern model are maintained. This construction only has an impact on the cross-covariances. 
Indeed, \cite{alegria2021bivariate}  showed that $\widetilde{K}_{12}(x)$ is not a monotonically decreasing function of $x$. More precisely, the cross-covariance can attain its maximum value at a strictly positive distance. This property was called cross-dimple in \cite{alegria2021bivariate}.
The parameter $\xi$ regulates the intensity of the cross-dimple. Clearly, the traditional bivariate Mat\'ern model is a limit case of this construction ($\xi\rightarrow \infty$). Closed-form expressions for (\ref{dimple}) are provided by \cite{alegria2021bivariate}. Moreover, for $\nu_{12} = 1/2+n$, $n\in\mathbb{N}$, $\widetilde{K}_{12}(x)$ can be expressed in terms of error functions and exponential functions.
\end{enumerate}
\item Directions, shapes and curves. The Mat{\'e}rn model is central to the study of directional processes. \citep{banerjee_directional}  formalize the notions of directional finite difference processes and directional derivative processes with special emphasis on the Mat{\'e}rn covariance function. They provide complete distribution theory results under the assumptions of a stationary Gaussian process (with Mat{\'e}rn covariance) model either for the data or for spatial random effects. \\
\cite{banerjee2006bayesian} introduced Bayesian wombling to measure gradients related to curves through wombling boun\-da\-ries. The smoothness properties of the Mat{\'e}rn model are proved to be successful within such a framework.
Modeling approaches to temporal gradients using the Mat{\'e}rn model have been proposed by \cite{quick2013modeling}.

\end{enumerate}



\end{document}